\documentclass[12pt]{article}
\usepackage[utf8]{inputenc}
\usepackage{authblk}
\usepackage{graphicx}
\usepackage{amssymb,amsmath}
\usepackage{mathtools}
\usepackage{textcomp}
\usepackage{multirow}
\usepackage[english]{babel}
\usepackage{color}
\usepackage{pgfplots}
\usepackage{tikz}
\usepackage{overpic}
\usepackage{caption}
\usepackage{subcaption}
\usepackage{breqn}
\usepackage{hyperref}
\usepackage{todonotes}
\usepackage[
top    = 2.0cm, bottom = 2.0cm, left   = 2.5cm, right  = 3.0cm]{geometry}

\numberwithin{equation}{section}

\title{A phase-field approach to model evaporation from porous media: Modeling and upscaling}
\author{Tufan Ghosh\footnote{Corresponding author} \\
Institute for Modelling Hydraulic and Environmental Systems, University of Stuttgart, Pfaffenwaldring 61, 70569 Stuttgart, Germany \\
\url{tufan.ghosh@iws.uni-stuttgart.de} \\
\and
Carina Bringedal \\
Department of Computer science, Electrical engineering and Mathematical sciences, Western Norway University of Applied Sciences, Inndalsveien 28, Bergen, 5063, Norway \\
\url{Carina.Bringedal@hvl.no} \\
\and
Christian Rohde \\
Institute of Applied Analysis and Numerical Simulation, University of Stuttgart, Pfaffenwaldring 57, 70569 Stuttgart, Germany \\
\url{christian.rohde@mathematik.uni-stuttgart.de}  \\
\and 
Rainer Helmig \\
Institute for Modelling Hydraulic and Environmental Systems, University of Stuttgart, Pfaffenwaldring 61, 70569 Stuttgart, Germany \\
\url{rainer.helmig@iws.uni-stuttgart.de} 
}

\newtheorem{remark}{Remark}[section]

\begin{document}

\maketitle

\begin{abstract}
\noindent We develop a phase-field model for evaporation from a porous medium by explicitly considering a vapor component together with the liquid and gas phases in the system. The phase-field model consists of the conservation of mass (for phases and vapor component), momentum, and energy. In addition, the evolution of the phase field is described by the Allen-Cahn equation. In the limit of vanishing interface width, matched asymptotic expansions reveal that the phase-field model reduces to the sharp-interface model with all the relevant transmission conditions on the moving interface. An energy estimate is derived, which suggests that for the diffusion-dominated regime, energy always decreases with time. However, this is not trivial in the case of other regimes. Through numerical examples, we analyze the efficiency of the developed phase-field formulation in modeling the evaporation process. We observe that our formulation is able to capture shrinking liquid droplet, in other words evaporation. Further, the phase-field model is upscaled to the Darcy scale using periodic homogenization for the diffusion-dominated regime. The effective parameters at the Darcy scale are connected to the pore scale through corresponding cell problems.
\end{abstract}

\section{Introduction}\label{Sec:1}
Evaporation in porous media is encountered in various industrial and environmental processes. As an example, in a geothermal reservoir evaporation due to cooling and low-pressure effect can lead to plant closing \cite{Bodvarsson1982}. On the other hand, drying in soils can result in salt precipitation which affects cultivation \cite{Ott2014}. Evaporation can be described as mass transfer from the liquid phase (water) to gas phase (air), along with corresponding momentum and energy transfer between the phases. One may note that the liquid-gas interface will also move, due to mass transfer. In this study, we aim to account for the evolving liquid-gas interface, when describing evaporation in a porous medium.

Modeling evaporation in porous media has been considered by many research groups, although mainly at the Darcy scale only (see Ref. \cite{Mosthaf2011} and references therein), with various heuristic approaches and strategies to determine the needed REV-scale parameters \cite{Ahrenholz2011,Mosthaf2014,Nuske2014}. A general approach for simulating non-isothermal multi-phase flow considering evaporation in porous media was developed in Refs. \cite{Class2002b} and \cite{Class2002a}. A simplified evaporation process is modeled in Refs. \cite{Class2002b} and \cite{Class2002a}, assuming local thermodynamic equilibrium in the sense that all the phases have the same temperature, which allowed them to consider only one energy balance for the entire porous medium. In the context of evaporation from a porous medium, Mosthaf \emph{et al.} \cite{Mosthaf2011} coupled two-phase porous-media flow with single-phase free flow. Later Mosthaf \emph{et al.} \cite{Mosthaf2014} extended the work in Ref. \cite{Mosthaf2011} by including a more refined evaporation rate by quantifying the vapor transfer between the porous medium and the free flow at Darcy scale. These models did include the vapor as a component of the gas, and accounted for non-ideal behavior of the gas phase. However, these approaches are not connected to the pore scale, rather postulated at the Darcy scale. Though reasonable and plausible, a detailed derivation from the pore-scale processes is certainly missing.

In a pore-scale formulation the in-situ phases are considered separately; \emph{i.e.,} the solid and saturating fluids are modeled through individual conservation equations and transfers between the phases across the evolving interfaces are explicitly accounted for. Pore-scale formulations with evolving interfaces between the in-situ phases have been widely used for \emph{e.g.} reactive transport where the fluid-solid interface evolves due to salt precipitation and dissolution \cite{vanNoorden2009b}. Usually, fluids are considered to be viscous and having low flow rates in porous media, thus the Navier-Stokes equations with low Reynolds number or the Stokes equations are suitable for fluid flow modeling at the pore scale \cite{Helmig1997}. Energy transfer in the phases takes place mainly through convection (fluids only) and conduction (fluids and solid). Across the interfaces, transmission conditions are prescribed to ensure conservation of mass, momentum and energy. Due to low flow rates, it is usually assumed that the phases share the same temperature on their common interface \cite{Bringedal2015}, which will be referred to as pore-scale thermal equilibrium. However, it may be noted that one can still have thermal non-equilibrium at the Darcy scale, since the averaged Darcy-scale temperature can still be different \cite{Bringedal2016a,Nuske2014}.

A large variety of models describing evaporation exist, depending on which processes are accounted for, and which flow regime and application is considered. Considering viscous fluids, Juric and Tryggvason \cite{Juric1998} formulated a single-field model for conservation of mass, momentum and energy including liquid-vapor phase change accounting for surface tension, interphase mass transfer and latent heat between the fluids. Gibou \emph{et al.} \cite{Gibou2007} considered an evaporation process incorporating jump conditions for mass, momentum and energy on the evolving interface. However, the models in Refs. \cite{Gibou2007} and \cite{Juric1998} considered incompressible fluids only, and that the gas phase consisted of only vapor.

How to describe evolving interfaces? This is one major challenge in formulating a pore-scale model involving evolving interfaces. One can easily differentiate between two obvious approaches; namely considering the interface as a sharp transition or as a diffuse transition between the phases. In case of sharp interfaces, the location of the interface can be followed either using the exact coordinates of the interface as unknowns, or implicitly through a level-set approach. When the evolving interface is considered to be sharp, the conservation equations for each phase are formulated in time-dependent domains. Hence, the conservation equations for each phase and the corresponding boundary conditions on the evolving interface are implicitly coupled with the equation describing the evolving interface. Refs. \cite{Bringedal2016b}, \cite{Gibou2007} and \cite{vanNoorden2009b} used level-set techniques to describe the sharp interface, whereas Bringedal \emph{et al.} \cite{Bringedal2015} could take the advantage of a simplified geometry and used the coordinates of the interface as the unknown. Note that the models described in the above studies have inherent discontinuities on the sharp interface, which makes analysis and implementation troublesome.

For diffuse-interface approaches, the transition between the two phases are considered smooth and of non-zero width $\lambda > 0$. Several approaches to mathematically describe the transition exist, and one such is applying a phase-field variable, which is a smooth phase indicator function approaching 1 in one phase, and 0 in the other phase. The evolving domains and interface are then handled by a phase-field equation; that is, an Allen-Cahn or Cahn-Hilliard equation. In a phase-field formulation, a common conservation equation for the two phases is formulated for each conserved variable (\emph{e.g.} energy), where the phase-field variable is incorporated into the equation to indicate which phase is considered. The transmission conditions between the phases are also included into the reformulated conservation equations. This way, the phase-field formulation is on a stationary domain and do not include discontinuities, which makes it more amendable for homogenization and numerical implementation. However, a phase-field formulation of a given sharp-interface problem is not unique, and it must be shown that the phase-field formulation captures the sharp-interface formulation when letting the width of the diffuse interface, $\lambda$, approach zero. The sharp-interface limit of a phase-field formulation can be found by applying matched asymptotic expansions \cite{Caginalp1988}.

Although phase-field formulations can make homogenization and numerical implementation easier, finding them is not straightforward. In a pore-scale setting, phase-field formulations with an evolving solid-fluid interface incorporating surface tension was developed by Redeker \emph{et al.} \cite{Redeker2016} and later combined with single-phase flow by Bringedal \emph{et al.} \cite{Bringedal2020}. For both these works, the sharp-interface limit was identified, and both models were also successfully upscaled to Darcy scale and implemented numerically. A thermodynamically consistent phase-field formulation for incompressible two-phase flow honoring the density difference between the fluids was formulated by Abels \emph{et al.} \cite{Abels2012}. They could also identify the sharp-interface limits for several versions of the phase-field equation. Phase-field formulations for two-phase flow with phase transitions and chemical reactions were developed in Refs. \cite{Aki2014}, \cite{Dreyer2014}, \cite{Witterstein2011} and \cite{Rohde2021}, where sharp-interface limits for certain regimes were found. However, all the above models are isothermal.

Berti and Giorgi \cite{Berti2009} considered a simplified, but thermodynamically consistent, phase-field formulation accounting for evaporation due to temperature and pressure changes. They accounted for the density changes in the vapor using a simplified equation of state. However, no sharp-interface limit was identified. Badillo \cite{Badillo2012} formulated a phase-field model for non-isothermal two-phase incompressible flow including evaporation, where surface-tension effects are taken into account for the evaporation rate. Using matched asymptotic analysis, Ref. \cite{Badillo2012} also showed how the formulation reduces to the corresponding sharp-interface model. Fabrizio \emph{et al.} \cite{Fabrizio2016} formulated a phase-field model for non-isothermal two-phase compressible flow where evaporation could take place. They accounted for the air component and varying densities due to evaporation, but did not analyze the sharp-interface limit of the model. Spinodal decomposition processes on the
pore scale are addressed for a non-isothermal fluid in Keim
\emph{et al.} \cite{Keim2023} using the Navier-Stokes Korteweg approach. Recently, Maier \emph{et al.} \cite{Maier2024} developed a Navier–Stokes Cahn–Hilliard model coupled with balance equations for heat and moisture to simulate the two-phase flow with phase change. They compared the phase-field model simulation with micro-fluidic experiments to analyze convective drying in porous media.

Pore-scale models are useful for understanding the underlying physics as they allow for a detailed description of the relevant processes. However, as the porous structure is highly detailed, direct numerical simulations of a domain interesting for real-life applications would be too computationally expensive. On the other hand, one is not necessarily interested in the detailed behavior, but mainly in the Darcy-scale effect of the pore-scale processes. In this case upscaling of the pore-scale formulation to Darcy scale can be applied. Several approaches for upscaling from pore to Darcy scale exist, and the most common ones are volume averaging (see Refs. \cite{Duval2004} and \cite{Niessner2009}), and homogenization relying on two-scale asymptotic expansions (see Ref. \cite{Hornung1996}), which will be considered in this study. We highlight that Duval \emph{et al.} \cite{Duval2004} formulated a pore-scale model for conservation of mass and energy on fluid domains that evolve due to evaporation. Using volume averaging they derived a Darcy-scale evaporation rate, which depends only on the temperature. Note that Ref. \cite{Duval2004} can be seen as a building block for adding processes and interactions in a pore-scale description. Later, through incorporating fluid flow and interfacial area as a variable in the volume averaging method, Niessner and Hassanizadeh \cite{Niessner2009} reformulated the model given in Ref. \cite{Duval2004}.

In case of applying homogenization through two-scale asymptotic expansions, a clear scale separation between pore-scale and Darcy-scale spatial variables is identified through the small number $\varepsilon$, which is the length scale ratio between the two scales. A local zoomed-in variable is applied for the pore scale, while an effective zoomed-out variable is used for the Darcy scale. The model unknowns are assumed to be expanded in terms of $\varepsilon$, using a series of functions depending on both spatial variables. For a periodic porous medium the effective equations rely on local cell problems in the pore-scale variable (see Ref. \cite{Hornung1996}).

For homogenization of two-phase flow in porous media, some choices must be made regarding the interaction between the flow and the evolving fluid-fluid interface to separate the scales. Possible approaches are to decouple the flow problem from the rest of the model by either applying a drift-approach as done in Ref. \cite{Schmuck2013}, or by taking advantage of the difference in time scales through scaling time with $\varepsilon$ as considered in Refs. \cite{Daly2015}, \cite{Metzger2021}, and \cite{Picchi2018} for different geometric settings and assumptions on flow regimes. In the latter approach, a Darcy-like law is obtained only when the fluid-fluid interface moves very slowly due to the flow.

Upscaling a sharp-interface formulation with an evolving boundary can be tedious as the evolving boundary needs an asymptotic expansion. In addition to this the time-dependent domains must also be accounted for when integrating in the pore scale. Although the transmission conditions on the evolving boundary need special attention, it has been successfully performed in various studies, see Refs. \cite{Bringedal2015,Bringedal2016b,Bringedal2016a,vanNoorden2009b} and \cite{vanNoorden2009a}.
Redeker \emph{et al.} \cite{Redeker2016} and Bringedal \emph{et al.} \cite{Bringedal2020}, showed that upscaling a phase-field formulation is more straightforward, since the domain of interest is stationary and the transmission conditions are already incorporated in the phase-field model equations. One may note that the phase field must also be expanded in terms of the scale separator $\varepsilon$, hence a choice must be made regarding how large the diffuse interface width is compared to the length scale at the pore scale. The choice made here will affect the resulting upscaled phase-field equation.

In this study, we develop a phase-field model for evaporation in porous media taking into account the necessary processes and interactions at the pore scale. Then we will apply periodic homogenization strategies to derive a Darcy-scale counterpart. The manuscript is structured as follows. In \hyperref[Sec:2]{Section 2} we describe the sharp-interface model considering the relevant processes at the pore scale. Then in \hyperref[Sec:3]{Section 3}, a phase-field formulation is introduced, based on the sharp-interface formulation. In addition to this an energy inequality for this model is developed. Next, in \hyperref[Sec:4]{Section 4} we derive the sharp-interface limit of the phase-field formulation, when the diffuse interface width approaches to zero. Further, the phase-field formulation developed in \hyperref[Sec:3]{Section 3} is solved numerically to analyze droplet evaporation from a rectangular channel in \hyperref[Sec:5]{Section 5}. Later, using periodic homogenization techniques a Darcy-scale counterpart of the phase-field model is developed in \hyperref[Sec:6]{Section 6} for the diffusion-dominated regime. Finally, we conclude this article with brief summary and outlook in \hyperref[Sec:7]{Section 7}.

\section{Sharp-Interface Formulation}\label{Sec:2}
\begin{figure}[h]
    \centering
    \includegraphics[scale=0.9]{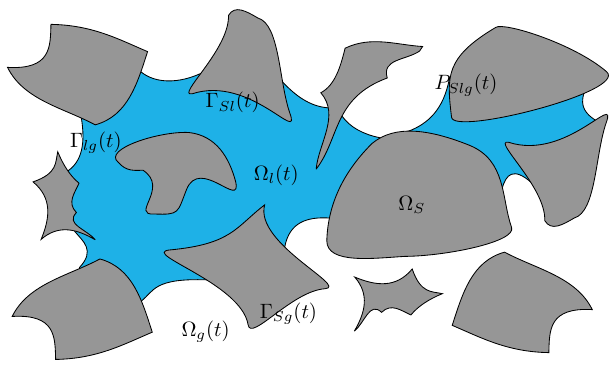}
    \caption{Schematic of a three-phase domain consisting of solid (gray), liquid (blue), and gas (white) phases.}
    \label{F1}
\end{figure}

In this section, we formulate the sharp-interface model with an evolving (free) boundary, which is going to be approximated by the proposed phase-field model in a later section. Here we consider a porous medium on the pore scale, saturated with two different fluids. Then, one has to take into account three phases namely, liquid, gas, and solid for evaporation from porous media. Here, the solid phase $\Omega_S$ is non-evolving and a stationary domain. The liquid and gas regions are evolving with time, and are denoted by $\Omega_l (t)$ and $\Omega_g (t)$, respectively. However, the combined liquid and gas phases' domain (i.e., the fluid domain) are stationary and denoted by $\Omega_F$. We denote the boundary between $\Omega_l (t)$ and $\Omega_g (t)$ with $\Gamma_{lg} (t)$. We may note that the interface between e.g. solid and liquid is evolving in the sense that the location of the liquid phase, and hence where it meets the solid, changes. The boundary between the solid and either the liquid or gas is denoted $\Gamma_{Sl} (t)$ and $\Gamma_{Sg} (t)$, respectively. In addition, the combined boundary between the solid and the fluids ($\Gamma_{Sl} (t)$ and $\Gamma_{Sg} (t)$) is denoted by $\Gamma_S$, which is stationary. We denote the triple-point where solid, liquid, and gas meet as $P_{Slg} (t)$. See \hyperref[F1]{Figure 1} for an overview.

In the following, we formulate the conservation of mass, momentum, and energy balances for the three phases, making relevant assumptions for the considered physical setting. In general, we assume that processes are slow in the sense that velocities are low and that local thermodynamic equilibrium on the various interfaces can be assumed.

\subsection{Conservation of Mass}\label{Sec:2-1}
Since the solid phase is immobile and also chemically non-reactive, we do not have to consider its mass balance. It is simply reduced to ${\bf v}_S = {\bf 0}$ in $\Omega_S$; \emph{i.e.,} there is no velocity in the solid. In general, the liquid and gas phases both consist of vapor and air as components. 
However, through evaporation, the vapor (water) becomes a component in the gas phase. Hence, we need to consider the mass conservation of the components as well as the phases. The mass balance for the two fluid phases is formulated as
\begin{equation}\label{ES1-1}
\partial_t \rho_{\alpha} + \nabla \cdot ( \rho_{\alpha} {\bf v}_{\alpha} ) = 0 \qquad \qquad \text{in } \Omega_{\alpha}(t), ~~ \alpha = l, g,
\end{equation}
\begin{equation}\label{ES1-2}
(\rho_l {\bf v}_l - \rho_g {\bf v}_g) \cdot {\bf n} = (\rho_l - \rho_g) v_n \qquad \qquad \text{on } \Gamma_{lg}(t),
\end{equation}
where $\rho_{\alpha}$ and ${\bf v}_{\alpha}$ are the density and velocity of the phase $\alpha$, respectively. The condition on $\Gamma_{lg} (t)$ is a jump condition ensuring the conservation of mass of the phases when there is a mass transfer (in the present study evaporation) between the phases. The normal velocity of the evolving fluid-fluid interface $v_n$ is a priori unknown and is defined with a positive direction into the gas. Note that a positive $v_n$ leads to condensation, while a negative $v_n$ corresponds to evaporation. Moreover, the normal vector ${\bf n}$ of the fluid-fluid interface is pointing into the gas phase. We define $\dot{m}$ to be the mass transfer rate between liquid and gas, as
\begin{equation}\label{ES1-3}
\dot{m} := \rho_{\alpha} ( {\bf v}_{\alpha} \cdot {\bf n} - v_n ),
\end{equation}
which is positive when the mass goes from liquid to gas.

Now for the mass conservation of the vapor component in the system, we consider the conservation of the moles, and hence formulate conservation using mole fractions (see Ref. \cite{Class2002a}). Diffusion inside the phases is taken into account in the governing equation using Fick's law of diffusion.
\begin{equation}\label{ES1-4}
\partial_t (\rho_{\alpha} \chi_{\alpha}^{v}) + \nabla \cdot (\rho_{\alpha} \chi_{\alpha}^{v} {\bf v}_{\alpha}) = \nabla \cdot ({\rm D}_{\alpha}^{v} \rho_{\alpha} \nabla \chi_{\alpha}^{v}) \qquad \text{in } \Omega_\alpha (t), ~~ \alpha = l, g
\end{equation}
\begin{equation}\label{ES1-5}
\left\{ \left( \rho_{l} \chi_l^{v} {\bf v}_l - {\rm D}_l^{v} \rho_{l} \nabla \chi_l^{v}  \right) - \left( \rho_{g} \chi_g^{v} {\bf v}_g - {\rm D}_g^{v} \rho_{g} \nabla \chi_g^{v}  \right) \right\} \cdot {\bf n} = (\rho_l \chi_l^{v} - \rho_{g} \chi_g^{v}) v_n \quad \text{on } \Gamma_{lg}(t),
\end{equation}
where $\chi_{\alpha}^{v}$ and ${\rm D}_{\alpha}^{v}$ are the mole fraction and the diffusivity of the vapor component in the phase $\alpha$.

\subsection{Conservation of Momentum}\label{Sec:2-2}
We assume that both the fluids (liquid and gas) are Newtonian fluids, and hence the conservation of momentum follows from the well-known Navier-Stokes equations. In addition, we also assume that the viscosity of the fluids is constant.
\begin{equation}\label{ES1-6}
\partial_t (\rho_{\alpha} {\bf v}_{\alpha}) + \nabla \cdot (\rho_{\alpha} {\bf v}_{\alpha} \otimes {\bf v}_{\alpha}) = - \nabla p_{\alpha} + \nabla \cdot \boldsymbol{\tau}_{\alpha} + \rho_{\alpha} {\bf g} \qquad \text{in } \Omega_{\alpha} (t), ~~ \alpha = l, g,
\end{equation}
\begin{equation}\label{ES1-7}
\left\{ - (p_l - p_g) {\rm {\bf I}} + ( \boldsymbol{\tau}_l - \boldsymbol{\tau}_g) \right\} \cdot {\bf n} = \dot{m}( {\rm {\bf v}}_l - {\rm {\bf v}}_g ) - \sigma \kappa {\bf n} \qquad \text{on } \Gamma_{lg} (t),
\end{equation}
\begin{equation}\label{ES1-8}
{\rm {\bf v}}_l \cdot {\rm {\bf t}}_i = {\rm {\bf v}}_g \cdot {\rm {\bf t}}_i \qquad \text{on } \Gamma_{lg} (t).
\end{equation}
Here $p_{\alpha}$ denotes the pressure of the phase $\alpha$, and
\begin{equation}\label{ES1-9}
\boldsymbol{\tau}_{\alpha} := \mu_{\alpha} \left( \nabla {\rm {\bf v}}_{\alpha} + \nabla {\rm {\bf v}}_{\alpha}^{T} \right) + \xi_{\alpha} (\nabla \cdot {\rm {\bf v}}_{\alpha}) {\rm {\bf I}},
\end{equation}
is the stress tensor, where $\mu_{\alpha}$ and $\xi_{\alpha}$ are the dynamic and bulk viscosities of the phase $\alpha$. Further, $\sigma$ is the surface tension between the fluids, and $\kappa = \nabla \cdot {\bf n}$ is the curvature of the fluids and is positive when the liquid forms droplets. In addition, we also assume there is no slip between the fluids, which is described by the relation given in \eqref{ES1-8}. Note that ${\bf t}_i$ is a tangential vector to the fluid-fluid interface (for a two-dimensional model we have only one such vector, while for a three-dimensional case, there are two orthogonal tangent vectors, $i = 1, 2$). Now for the solid phase, the conservation of momentum simply states
\begin{equation}\label{ES1-10}
{\bf v}_S = {\bf 0} \qquad \qquad \text{in } \Omega_S.
\end{equation}
Thus when we apply a no-slip boundary condition between the fluids and the solid, it becomes
\begin{equation}\label{ES1-11}
{\bf v}_{\alpha} = {\bf 0} \qquad \qquad \text{on } \Gamma_{S\alpha} (t).
\end{equation}
Note that no-slip boundary conditions for the sharp-interface formulation of two-phase flow are singular at the fluid-fluid-solid contact line, and usually require a form of regularized slip condition for the contact line $\Gamma_{S\alpha}$. However, since this will not be an issue for the phase-field formulation, which is the main focus of this work, we do not formulate such regularized slip conditions here but refer to e.g. Ref. \cite{Lunowa2021} for such an approach. For the sake of simplicity, here we assume a (neutral) $90^{\rm o}$ contact angle between the fluids and the solid.

\subsection{Conservation of Energy}\label{Sec:2-3}
The general form for conservation of energy in the fluid phase $\alpha$ is given by
\begin{equation}\label{ES1-12}
\begin{aligned}
\partial_t \left( \rho_{\alpha} \left( u_{\alpha} + \frac{v^2_{\alpha}}{2} \right) \right) + \nabla \cdot \left( \rho_{\alpha} \left( u_{\alpha} + \frac{v^2_{\alpha}}{2} \right) {\bf v}_{\alpha} \right) + \nabla \cdot ( p_{\alpha} {\bf v}_{\alpha} ) & \\
= \nabla \cdot ( k_{\alpha} \nabla {\rm T}_{\alpha} ) + \nabla \cdot ({\bf v}_{\alpha} \cdot \boldsymbol{\tau}_{\alpha}) + \rho_{\alpha} {\rm {\bf v}}_{\alpha} \cdot {\rm {\bf g}} & \quad \text{in } \Omega_{\alpha}(t).
\end{aligned}
\end{equation}
Here, $u_{\alpha}$, $v_{\alpha}^{2} = \|\mathbf v_\alpha\|^2$, $k_{\alpha}$ and ${\rm T}_{\alpha}$ are the internal energy per unit mass, the square of the velocity magnitude, the heat conductivity and the temperature of the phase $\alpha$, respectively. Note that the first term on the left-hand side of equation \eqref{ES1-12} represents energy accumulation, the next two terms describe energy transport due to advection and volume change, respectively; while the terms on the right-hand side represent energy transport due to conduction, viscous stress, and body forces respectively. We assume that all the phases follow Fourier's law for heat conduction.

One may note that equation \eqref{ES1-12} contains the equation for the conservation of kinetic energy. This can be easily derived from the momentum equation \eqref{ES1-6}, by taking the dot product of the velocity ${\rm {\bf v}}_{\alpha}$ with the equation \eqref{ES1-6}. Using some simple algebraic manipulations lead to
\begin{equation}\label{ES1-13}
\partial_t \left( \rho_{\alpha} \frac{ v_{\alpha}^{2}}{2} \right) + \nabla \cdot \left( \rho_{\alpha} {\rm {\bf v}}_{\alpha} \frac{v_{\alpha}^{2}}{2} \right) = - {\rm {\bf v}}_{\alpha} \cdot \nabla p_{\alpha} + {\rm {\bf v}}_{\alpha} \cdot \left( \nabla \cdot \boldsymbol{\tau} \right) + \rho_{\alpha} {\rm {\bf v}}_{\alpha} \cdot {\bf g} \quad \text{in } \Omega_{\alpha}(t).
\end{equation}
We note that all of the work of the body forces in the total energy equation \eqref{ES1-12} goes to change the kinetic energy. Then, the rest of the right-hand side of equation \eqref{ES1-13} consists of the parts of the surface force work term as affecting the kinetic energy. Now, subtracting equation \eqref{ES1-13} from equation \eqref{ES1-12} gives the internal energy equation
\begin{equation}\label{ES1-14}
\partial_t (\rho_{\alpha} u_{\alpha}) + \nabla \cdot (\rho_{\alpha} u_{\alpha} {\bf v}_{\alpha}) = \nabla \cdot ( k_{\alpha} \nabla {\rm T}_{\alpha}) - p_{\alpha} \nabla \cdot {\bf v}_{\alpha} + \boldsymbol{\tau}_{\alpha} : \nabla {\bf v}_{\alpha}  \quad \text{in } \Omega_{\alpha}(t).
\end{equation}
Let $h$ be the specific enthalpy defined as $h := u + \frac{p}{\rho}$, then the above equation (\ref{ES1-14}) can be rewritten as
\begin{equation}\label{ES1-15}
\partial_t (\rho_{\alpha} u_{\alpha}) + \nabla \cdot (\rho_{\alpha} h_{\alpha} {\bf v}_{\alpha}) = \nabla \cdot ( k_{\alpha} \nabla {\rm T}_{\alpha}) + {\bf v}_{\alpha} \cdot \nabla p_{\alpha} + \boldsymbol{\tau}_{\alpha} : \nabla {\bf v}_{\alpha} \quad \text{in } \Omega_{\alpha}(t).
\end{equation}
One may note that the above equation \eqref{ES1-15} is supplemented by the equations of state relating the internal energy with the temperature for each of the phases (see Ref.~\cite{Batzle1992}). We assume local thermal equilibrium on the fluid-fluid interface which is expressed as
\begin{equation}\label{ES1-16}
{\rm T}_l = {\rm T}_g \qquad \qquad \text{on } \Gamma_{lg} (t).
\end{equation}
We can formulate the jump condition on the fluid-fluid interface by
\begin{equation}\label{ES1-17}
\left\{ - \left( k_l \nabla {\rm T}_l - k_g \nabla {\rm T}_g \right) + \left( \rho_l u_l {\bf v}_l - \rho_g u_g {\bf v}_g \right) \right\} \cdot {\bf n} = \left( \rho_l u_l - \rho_g u_g \right) v_n \qquad \text{on } \Gamma_{lg} (t).
\end{equation}
Using the definition of $\dot{m}$, the mass transfer rate between the liquid and gas, the last relation (\ref{ES1-17}) can be rewritten as
\begin{equation}\label{ES1-18}
\left( k_l \nabla {\rm T}_l - k_g \nabla {\rm T}_g \right) \cdot {\bf n} = \dot{m} \left( u_l - u_g \right) \qquad \text{on } \Gamma_{lg} (t).
\end{equation}
Note that the above relation (\ref{ES1-18}) is very often approximated to
\begin{equation}\label{ES1-19}
\left( k_l \nabla {\rm T}_l - k_g \nabla {\rm T}_g  \right) \cdot {\bf n} \approx \dot{m} \left( h_l - h_g \right) \qquad \text{on } \Gamma_{lg} (t).
\end{equation}
The latent heat of evaporation $\mathcal{L}$; \emph{i.e.} the amount of enthalpy that is consumed when evaporation occurs, is defined by
\begin{equation}\label{ES1-20}
\mathcal{L} := h_{g} - h_{l},
\end{equation}
where enthalpy values are taken from the fluid-fluid interface. Usually one can use the saturated values for enthalpy on the interface. Then the latent heat of evaporation is expressed as
\begin{equation}\label{ES1-21}
\mathcal{L} = h_{g}^{sat} - h_{l}^{sat}.
\end{equation}
Hence, the relation (\ref{ES1-19}) can be rewritten as
\begin{equation}\label{ES1-22}
\left( k_g \nabla {\rm T}_g - k_l \nabla {\rm T}_l  \right) \cdot {\bf n} \approx \dot{m} \mathcal{L} \qquad \text{on } \Gamma_{lg} (t).
\end{equation}
Furthermore, the internal energy of the solid phase depends only on the temperature, hence $u_S = c_{p, S} {\rm T}_S$, and the corresponding energy balance is given by
\begin{equation}\label{ES1-23}
\partial_t (\rho_S c_{p,S} {\rm T}_S) = \nabla \cdot ( k_S \nabla {\rm T}_S) \qquad \text{in } \Omega_S.
\end{equation}
Here, $\rho_S$ and $c_{p, S}$ are the density and heat capacity of the solid phase, respectively. In the present study, for simplicity, we choose both the density and heat capacity of the solid as constant. We assume continuity of temperature and heat flux on the fluid-solid interfaces, \emph{i.e.,}
\begin{equation}\label{ES1-24}
{\rm T}_S = {\rm T}_{\alpha} \qquad \qquad \text{on } \Gamma_{S\alpha} (t),
\end{equation}
\begin{equation}\label{ES1-25}
k_S \nabla {\rm T}_S \cdot {\bf n}_S = k_{\alpha} \nabla {\rm T}_{\alpha} \cdot {\bf n}_S \qquad \text{on } \Gamma_{S\alpha} (t).
\end{equation}
Note that the convective part of the fluid heat flux disappears in \eqref{ES1-25}, since ${\bf v}_{\alpha} = 0$ on $\Gamma_{S\alpha} (t)$ (see relation \eqref{ES1-11}).

\subsection{Evaporation Rate}\label{Sec:2-4}
How does the normal velocity of the interface $v_n$ evolve? 
In connection with this, the only information we have is that this velocity $v_n$ is directly related to fluid flow and mass transfer across the interface. A usual assumption in modeling the evaporation process is that the evaporation is an equilibrium mechanism, \emph{i.e.,} the gas phase is always fully saturated with vapor near the fluid-fluid interface. In this case, one can easily calculate how much water needs to evaporate to (re)gain fully saturated conditions. However, this suggests that the normal velocity would be infinite, thus for our modeling perspective it is more beneficial to rather formulate a large evaporation rate, which would indicate that fully saturated conditions are achieved fast but not instantaneously. We formulate this through an evaporation rate $f (\chi^v_l, T_l)$, as follows
\begin{equation}\label{ES1-26}
\frac{\left( \rho_l + \rho_g \right)}{2} \left( v_{n} - {\bf v}_{l} \cdot {\bf n} \right) = - f(\chi_l^v, {\rm T}_l) \quad \text{on } \Gamma_{lg} (t).
\end{equation}

The kinetic theory of gases states that the rate of evaporation is directly linked to the thermodynamic non-equilibrium between the saturated vapor pressure $p_{sat}$ and the actual vapor pressure $p_v$ on the liquid-gas interface. Using the Hertz-Knudsen equation \cite{Knudsen1934, Doss2023}, one can estimate the local values of $f$ as follows
\begin{equation}\label{ES1-27}
f(\chi_l^v, {\rm T}_l) = \sqrt{\frac{M_l}{2 \pi \mathcal{R} {\rm T}_l}} \left( p_{sat} - p_v \right),
\end{equation}
where $M_l$ and $\mathcal{R}$ are the molar mass of liquid and the ideal gas constant, respectively. Then, we use Raoult's law to estimate saturated vapor pressure $p_{sat}$ in terms of the saturated vapor pressure of pure liquid $p_{l}^{sat}$, as $p_{sat} = \chi_l^v p_{l}^{sat}$. Further, one can use the August-Roche-Magnus equation \cite{Alduchov1996} to compute $p_{l}^{sat}$ in Pascal $(Pa)$, given by
\begin{equation}\label{ES1-28}
    p_{l}^{sat} = 610.94 \exp{\left( \frac{17.625 \Theta_l}{\Theta_l + 243.04} \right)}, 
\end{equation}
where $\Theta_l = {\rm T}_l - 273.15 K$ is the local liquid temperature in Celsius. Moreover, one can compute actual vapor pressure using ideal gas law $p_v M_l = \rho_v \mathcal{R} {\rm T}_l$.

\section{Phase-field Formulation}\label{Sec:3}
\begin{figure}[h]
  \centering
  \begin{overpic}[scale = 0.5]{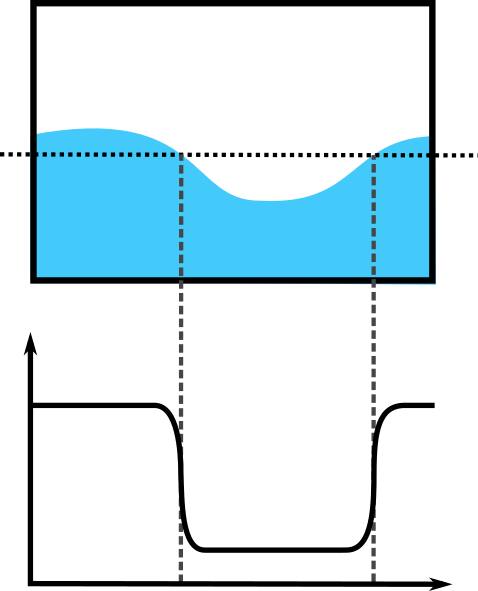}
      \put(40,82){\color{black}{\it Gas}}
      \put(35,58){\color{black}{\it Liquid}}
      \put(24,38){\color{black}{$\mathcal{O}(\lambda)$}}
      \put(25,32){$\boldsymbol{\longleftrightarrow}$}
      \put(0,29){\color{black}{1}}
      \put(0,6){\color{black}{0}}
  \end{overpic}
      \caption{A schematic representation of the phase field.}
      \label{F2}
\end{figure}

A useful alternative to the sharp-interface formulation described in the previous section is to consider a phase-field formulation. Note that the phase field is an approximate phase indicator function, with a smooth transition from one phase to the other. As a consequence, it creates a diffuse interface layer of finite width. In the present study, we adopt a phase-field formulation to describe the governing equations in the fluid domain $\Omega_F$ as a whole (\emph{i.e.,} liquid and gas together). To be specific, the dimensionless phase field $\phi$ approaches 0 in the gas phase, and 1 in the liquid phase, and has a smooth transition of (dimensional) width $\mathcal{O}(\lambda)$ separating the phases, as seen in \hyperref[F2]{Figure 2}. Here, $\lambda > 0$ indicates the thickness of the diffuse transition layer and is hence a phase-field parameter and does not represent anything physical. We will be required to arrive at the sharp-interface model in the limit of $\lambda \rightarrow 0$. Hence, the sharp interface between the two phases is approximated by a smooth, diffuse transition region in the phase-field formulation. The main advantage of this phase-field formulation is that we are now dealing with the stationary fluid domain $\Omega_F$, and no longer with the time-dependent (evolving) domains ($\Omega_l (t)$ and $\Omega_g (t)$). In the following, we present the phase-field model for evaporation in porous media.

\subsection{Phase-field Equation}\label{Sec:3-1}
We define the evolution of the phase field $\phi$ through an Allen-Cahn type equation given below:
\begin{dmath}\label{ES2-1}
\rho \left( \partial_{t} \phi + \nabla \cdot \left( \phi {\rm {\bf v}} \right) \right) = \frac{\gamma}{\nu} \left\{ \nabla^2 \phi - \lambda^{-2} P'(\phi) \right\} - \sqrt{2} \lambda^{-1} \phi (1 - \phi) f(\chi^{v}, {\rm T}) \qquad \text{in } \Omega_F,
\end{dmath}
where, $\rho = \phi \rho_l + (1 - \phi) \rho_g$ and ${\rm {\bf v}}$ are the density and velocity of the composite phase, respectively; $\gamma$ and $\nu$ are parameters having dimensions ${\rm kg ~ s^{-2}}$ and ${\rm m ~ s^{-1}}$, respectively. Here, $P(\phi) = \phi^2 (1 - \phi)^2$ is the double-well potential.  It may be observed that $P'(\phi)$ dominates the equation \eqref{ES2-3} if $\lambda$ is small, implying that $\phi$ approaches one of the three equilibrium values 0, $\frac{1}{2}$, 1. We later show that $\phi = \frac{1}{2}$ is an unstable equilibrium. Hence, the double-well potential $P(\phi)$ ensures that the phase field attains values (close to) 0 and 1 for $\lambda \rightarrow 0$. Here, the evaporation rate $f(\chi^{v}, {\rm T})$ is similar to one defined in the sharp-interface formulation. However, here it is valid in the entire domain instead of the local behavior in the sharp-interface formulation. Note that the factor $\phi (1 - \phi)$ ensures that the evaporation rate is non-zero only in the diffusive transition layer.

One can incorporate energetically consistent static contact angles through $\sigma \nabla \phi \cdot {\bf n} = - g'(\phi)$, where $g(\phi)$ is the wall free energy, which is connected to the microscopic contact angle at the wall \cite{Jacqmin2000}. Note that $g(\phi) = 0$ translates to a $90^{\rm o}$ (neutral) contact angle. Hence, in this case, we have the homogeneous Neumann condition on the fluid-solid interface, {\it i.e.,}
\begin{equation}\label{ES2-2}
\nabla \phi \cdot {\bf n} = 0 \qquad \text{on } \Gamma_S.
\end{equation}

The effect of the non-neutral contact angle is not the focus of this study. We would like to refer interested readers to \cite{Kelm2024} for more details on how to include non-neutral contact angles explicitly in the phase-field formulation and the corresponding upscaling procedure through appropriate scaling. 


\subsection{Mass Conservation Equations}\label{Sec:3-2}
Taking into account the assumptions made on \hyperref[Sec:2]{Section 2}, the mixture mass conservation equation for liquid and gas together along with the vapor component balance for the entire fluid system are expressed as follows
\begin{equation}\label{ES2-3}
\partial_{t} \rho + \nabla \cdot \left( \rho {\rm {\bf v}} \right) = 0 \qquad \text{in } \Omega_F,
\end{equation}
\begin{equation}\label{ES2-4}
\partial_{t} \left( \rho \chi^{v} \right) + \nabla \cdot \left( \rho \chi^{v} {\rm {\bf v}} \right) = \nabla \cdot \left( {\rm D}^{v} \rho \nabla \chi^{v} \right) \quad \text{in } \Omega_F,
\end{equation}
where, $\chi^v$ and ${\rm D}^v = \phi {\rm D}_l^v + (1 - \phi) {\rm D}_g^v$ denotes respectively the vapor mole fraction and the vapor diffusivity. Since there is no mass exchange between the solid and the fluids, thus we have a no-flux condition on the solid interfaces, which reads
\begin{equation}\label{ES2-5}
{\rm D}^{v} \rho \nabla \chi^{v} \cdot {\bf n} = 0 \qquad \text{on } \Gamma_S.
\end{equation}

\subsection{Momentum Balance Equation}\label{Sec:3-3}
Now, as we are considering both phases together as a mixture, one has to take into account the surface tension force between the liquid and the gas in the force balance. Therefore, the classical Navier-Stokes's equation is modified with the surface tension force term for the mixture momentum balance as
\begin{equation}\label{ES2-6}
\partial_{t} \left( \rho  {\rm {\bf v}} \right) +  \nabla \cdot \left( \rho {\rm {\bf v}} \otimes {\rm {\bf v}} \right) = - \nabla p + \nabla \cdot \boldsymbol{\tau} - \nabla \cdot \left( \lambda \sigma \nabla \phi \otimes \nabla \phi \right) + \rho {\rm {\bf g}} \qquad \text{in } \Omega_F,
\end{equation}
where, $p$, $\boldsymbol{\tau} = \mu \left( \nabla {\rm {\bf v}} + \nabla {\rm {\bf v}}^{T} \right) + \xi (\nabla \cdot {\rm {\bf v}}) {\rm {\bf I}}$, $\mu = \phi \mu_l + (1 - \phi) \mu_g$ and $\xi = \phi \xi_l + (1 - \phi) \xi_g$ are the pressure, the viscous stress tensor, dynamic and bulk viscosities of the mixture, respectively. Note that the solid phase has no velocity. Hence, the no-slip boundary condition on the solid interfaces reads
\begin{equation}\label{ES2-7}
{\rm {\bf v}} = {\bf 0} \qquad \text{on } \Gamma_S.
\end{equation}
One can certainly incorporate a Navier-slip condition instead of \eqref{ES2-7} through ${\rm {\bf v}} = - \ell_s (\nabla {\rm {\bf v}}_{\bf t} \cdot {\bf n}) {\bf t}$, where $\ell_s$ is the slip length, refer \cite{Kelm2024}.

\subsection{Energy Balance Equation}\label{Sec:3-4}
Similar to the momentum balance equation, here we will write the energy balance equation for the mixture. Thus we need to take into account the energy due to surface tension between the fluids. In doing so we modify the energy balance equation as follows.
\begin{equation}\label{ES2-8}
\partial_{t} \left( \rho u \right) + \nabla \cdot \left( \rho h {\rm {\bf v}} \right) = \nabla \cdot \left( k \nabla {\rm T} \right) + {\rm {\bf v}} \cdot \nabla p + \boldsymbol{\tau} : \nabla {\rm {\bf v}} + {\rm {\bf v}} \cdot \left\{ \nabla \cdot \left( \lambda \sigma \nabla \phi \otimes \nabla \phi \right) \right\} \quad \text{in } \Omega_F,
\end{equation}
where, $u = \phi u_l + (1 - \phi) u_g$ and $h = \phi h_l + (1 - \phi) h_g$ are the specific internal energy and the specific enthalpy of the mixture, respectively; $k = \phi k_l + (1 - \phi) k_g$ denotes the thermal conductivity of the mixture, ${\rm T}$ is the temperature of the mixture. One can note that the last term on the right-hand side of the above equation \eqref{ES2-8} is balancing the energy due to the surface tension force between the fluids. We assume continuity of temperature and heat flux between the mixture and the solid on the combined fluid-solid interface, \emph{i.e.,}
\begin{equation}\label{ES2-9}
{\rm T} = {\rm T}_S \qquad \qquad \text{on } \Gamma_{S},
\end{equation}
\begin{equation}\label{ES2-10}
k \nabla {\rm T} \cdot {\bf n} = k_S \nabla {\rm T}_S \cdot {\bf n} \qquad \text{on } \Gamma_{S}.
\end{equation}

\subsection{Free Energy Functional}\label{Sec:3-5}
The energy associated with the above phase-field model is given by
\begin{equation*}\label{ES2-11}
E = \int_{\Omega_F} \left( \frac{1}{2} \rho {\rm {\bf v}}^{2} + \gamma \lambda^{-1} P(\phi) + \frac{1}{2} \gamma \lambda \mid \nabla \phi \mid ^{2} + \rho u + \rho F(\rho, \phi) \right) d{\rm {\bf x}},
\end{equation*}
which is the sum of the kinetic energy, the free energy of the phase field, the internal energy, and the density energy. Here the density energy $\rho F(\rho, \phi)$ is defined as follows
\begin{equation*}\label{ES2-12}
\partial_{\phi} (\rho F(\rho, \phi)) = \sqrt{2} \nu \phi (1 - \phi) f(\chi^{v}, {\rm T}).
\end{equation*}
Straightforwardly one can compute
\begin{equation*}\label{ES2-13}
\begin{aligned}
\frac{d}{dt} E(t) = & ~ \frac{d}{dt} \int_{\Omega_F} \left( \frac{1}{2} \rho {\rm {\bf v}}^{2} + \gamma \lambda^{-1} P(\phi) + \frac{1}{2} \gamma \lambda \mid \nabla \phi \mid ^{2} + \rho u + \rho F(\rho, \phi) \right) d{\rm {\bf x}} \\
 = & \int_{\Omega_F} \left( \partial_t \left( \frac{1}{2} \rho {\rm {\bf v}}^{2} \right) + \partial_t \left( \gamma \lambda^{-1} P(\phi) + \frac{1}{2} \gamma \lambda \mid \nabla \phi \mid ^{2} + \rho F(\rho, \phi) \right) + \partial_t \left( \rho u \right) \right) d{\rm {\bf x}} \\
 = & \int_{\Omega_F} \left( {\rm {\bf v}} \cdot \partial_t (\rho {\rm {\bf v}}) - \frac{1}{2} {\rm {\bf v}}^2 \partial_t \rho + \left( \gamma \lambda^{-1} P'(\phi) - \gamma \lambda \nabla^{2} \phi + \partial_{\phi} (\rho F(\rho, \phi)) \right) \partial_t \phi \right) d{\rm {\bf x}} \\
 & + \int_{\Omega_F} \left( \partial_t (\rho u) + \partial_{\rho} (\rho F(\rho, \phi)) \partial_t \rho \right)  d{\rm {\bf x}}.
\end{aligned}
\end{equation*}
Now we utilize equations \eqref{ES2-1}, \eqref{ES2-3}, \eqref{ES2-4}, \eqref{ES2-6}, and \eqref{ES2-8} to substitute the time derivatives in the above volume integrals. Then using integration by parts and the divergence theorem (wherever applicable), along  with the no-slip boundary condition for ${\rm {\bf v}}$ and zero flux conditions for $\phi$ and ${\rm T}$ on a bounded domain $\Omega$, leads to
\begin{equation*}\label{ES2-14}
\frac{d}{dt} E(t) = \int_{\Omega_F} \left( {\rm {\bf v}} \cdot \rho {\rm {\bf g}} - \eta \nabla \cdot (\phi {\rm {\bf v}}) - \frac{\eta^{2}}{\rho \lambda \nu} - \partial_{\rho} \left( \rho F(\rho, \phi) \right) \nabla \cdot (\rho {\rm {\bf v}}) \right) d {\rm {\bf x}},
\end{equation*}
where $\eta = \gamma \lambda^{-1} P'(\phi) - \gamma \lambda \nabla^{2} \phi + \sqrt{2} \nu \phi (1 - \phi) f(p, T, \chi^{v})$. The first term on the right-hand side describes energy due to gravity forces, which would be balanced by potential energy. The second term and the fourth term can be negative or positive and thus potentially lead to increasing energy. Only the third term is always negative. Therefore, generally, the energy $E$ is not always decreasing over time. However, one may note that for zero or negligible velocity (\emph{i.e.,} for the diffusion-dominated regime), the energy $E$ decreases with time when the third term dominates.

\section{Sharp-Interface Limit}\label{Sec:4}
The phase-field model is defined in the entire (stationary) fluid domain and the free (evolving) boundary is approximated by a diffuse interface layer. Hence, this can be thought of as an approximation of the sharp-interface model. To show this, we here evaluate the limit of the phase-field model (see \hyperref[Sec:3]{Section 3}) as the width of the diffuse transition layer $\lambda$ approaches zero. We show that this limit is the sharp-interface model described in \hyperref[Sec:2]{Section 2}. To do this, we introduce a dimensionless parameter $\zeta = \frac{\lambda}{L}$ related to thickness of the diffuse interface layer, where $L$ is a typical length in the model. In the following, we investigate the behavior of the solution as $\zeta \rightarrow 0$. To do this, we express the unknowns as a series in terms of $\zeta$ and equate similar order terms. Following Ref. \cite{Caginalp1988}, we examine the solution behavior very close to the interface and far away from it by employing inner expansion very close to the diffuse interface, and outer expansion far away from it. In addition to these, we apply matching conditions to connect the inner and outer expansions.

Away from the interface, we consider the outer expansion of the unknowns. For $\phi$ this reads
\begin{equation}\label{ES3-1}
\phi^{\rm out}(t,{\bf x}) = \phi_{0}^{\rm out}(t,{\bf x}) + \zeta \phi_{1}^{\rm out}(t,{\bf x}) + \zeta^2 \phi_{2}^{\rm out}(t,{\bf x}) + \mathcal{O}(\zeta^3),
\end{equation}
and similarly for the other unknowns.

\begin{figure}[h!]
    \centering
    \begin{overpic}[width=0.35\textwidth]{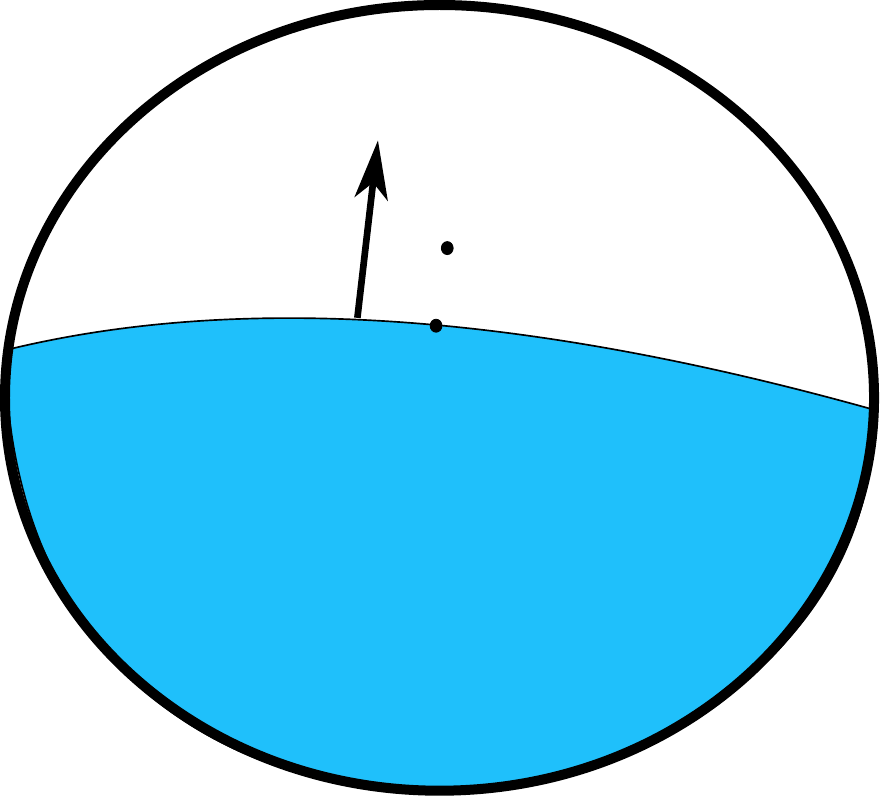}
    \put(52,46){\color{black} $\mathbf y_{\zeta}(t,\mathbf s)$}
    \put(54,61){\color{black} $\mathbf y_{\zeta}+r\mathbf n_{\zeta}(t,\mathbf s)$}
    \put(28,60){\color{black} $\mathbf n_{\zeta}$}
    \put(10,90){\color{black} $z=\zeta^{-1}r$}
    \put(30,25){\color{black} Liquid}
    \put(60,76){\color{black} Gas}
    \end{overpic}
    \caption{Local curvilinear coordinates near the interface.}
    \label{F3}
\end{figure}

Following Refs. \cite{Bringedal2020} and \cite{Redeker2016}, we switch to local coordinates near the diffuse interface where the inner expansion is valid. To be specific, we assume that $\Gamma_\zeta(t)$ denotes the set of points ${\bf y}_{\zeta} \in \Omega$ such that $\phi({\bf y}_{\zeta}, t) = 1/2$. These points vary with time $t$, and the dimensionless parameter $\zeta$ as the model depends on $\lambda = L \zeta$. Denoting ${\bf s}$ as the parameterization along $\Gamma_{\zeta}(t)$ and ${\bf n}_{\zeta}$ as the normal vector at $\Gamma_{\zeta}(t)$ pointing into the gas phase, we can define $r$ to be the signed distance to a point ${\bf x}$ near $\Gamma_{\zeta}(t)$. Hence, $r$ varies with ${\bf x}$ and $t$, and is considered positive in the gas phase. Thus we have
\begin{equation}\label{ES3-2}
{\bf x} = {\bf y}_{\zeta}(t, {\bf s}) + r {\bf n}_{\zeta}(t, {\bf s}),
\end{equation}
see \hyperref[F3]{Figure 3}. One can then show that
\[
\mid \nabla r \mid = 1, \quad \nabla r \cdot \nabla s_i = 0, \quad \partial_t r = - v_n, \quad \nabla^2 r = \frac{\kappa + 2 \Pi r}{1 + \kappa r + \Pi r^2},
\]
where $\kappa$ and $\Pi$ are the mean and Gaussian curvature of the moving interface, see Ref. \cite{Caginalp1988}. Next, we assume that ${\bf y}_{\zeta} = {\bf y}_{0} + \zeta {\bf y}_{1} + \cdots$, where ${\bf y}_{0}$ is a point on the interface $\Gamma_0^{\rm out} (t)$ defined through $\phi_0^{\rm out} = 1/2$. Analogously, we assume ${\bf n}_{\zeta} = {\bf n}_{0} + \mathcal{O}(\zeta)$, where ${\bf n}_{0}$ is the normal vector of $\Gamma_0^{\rm out} (t)$.

Now, we consider the inner expansions in terms of $z$ and ${\bf s}$, using the relation $z = \zeta^{-1}r$. This reads
\begin{equation}\label{ES3-3}
\phi^{\rm in}(t,{\bf x}) = \phi_{0}^{\rm in}(t,z,{\bf s}) + \zeta \phi_{1}^{\rm in}(t,z,{\bf s}) + \zeta^2 \phi_{2}^{\rm in}(t,z,{\bf s}) + \mathcal{O}(\zeta^3),
\end{equation}
and correspondingly for the other unknowns. Considering the scaling of the $z$ variable and the curvilinear coordinates defined in \eqref{ES3-2}, the derivatives can be rewritten as follows. For a generic variable $a$ or a generic vector-valued function ${\bf a}$, we obtain
\begin{equation}\label{ES3-4}
\partial_t a = - \zeta^{-1} v_{n,0} \partial_z a^{\rm in} + \left( \partial_t + \partial_t s \cdot \nabla_s \right) a^{\rm in} + \mathcal{O}(\zeta),
\end{equation}
\begin{equation}\label{ES3-5}
\nabla a = \zeta^{-1} \partial_z a^{\rm in} {\bf n}_0 + \nabla_{\Gamma} a^{\rm in} + \mathcal{O}(\zeta),
\end{equation}
\begin{equation}\label{ES3-6}
\nabla \cdot {\bf a} = \zeta^{-1} \partial_z {\bf a}^{\rm in} \cdot {\bf n}_0 + \nabla_{\Gamma} \cdot {\bf a}^{\rm in} + \mathcal{O}(\zeta),
\end{equation}
\begin{equation}\label{ES3-7}
\nabla {\bf a} = \zeta^{-1} \partial_z {\bf a}^{\rm in} \otimes {\bf n}_0 + \nabla_{\Gamma} {\bf a}^{\rm in} + \mathcal{O}(\zeta),
\end{equation}
\begin{equation}\label{ES3-8}
\nabla^2 a = \zeta^{-2} \partial_{zz} a^{\rm in} + \zeta^{-1} \kappa_0 \partial_z a^{\rm in} + \mathcal{O}(1),
\end{equation}
where $\kappa_0$ and $v_{n,0}$ are the mean curvature and normal velocity of the interface $\Gamma_0^{\rm out} (t)$. Note that, we have utilized the relations $\nabla^2 r = \kappa_0 + \mathcal{O}(\zeta)$ and $v_n = v_{n,0} + \mathcal{O}(\zeta)$ in order to re-write the above derivatives. Additionally, the properties $|\nabla r| = 1$ and $\nabla r \cdot \nabla s_i = 0$ are used to arrive at the last equality.

For the outer expansion and fixed $t$ and ${\bf s}$, we assume that ${\bf y}_{1/2 \pm}$ denote the limit $r \searrow 0$ (\emph{i.e.} from the gas phase), respectively $r \nearrow 0$ (from the liquid phase) of ${\bf x}$ rewritten in terms of the local coordinates \eqref{ES3-2}. We match the corresponding limit values of the outer and the inner expansions for $z \rightarrow \pm \infty$. Following Ref. \cite{Caginalp1988}, we explicitly assume that the phase field $\phi$ (and similarly for all other unknowns) satisfies the following matching conditions
\begin{subequations}\label{ES3-9}
\begin{equation}\label{ES3-9a}
\lim_{z \rightarrow \pm \infty} \phi_{0}^{\rm in} (t,z,{\bf s}) = \phi_{0}^{\rm out} (t,{\bf y}_{1/2\pm}),
\end{equation}
\begin{equation}\label{ES3-9b}
\lim_{z \rightarrow \pm \infty} \partial_{z} \phi_{0}^{\rm in} (t,z,{\bf s}) = 0,
\end{equation}
\begin{equation}\label{ES3-9c}
\lim_{z \rightarrow \pm \infty} \partial_{z} \phi_{1}^{\rm in} (t,z,{\bf s}) = \nabla \phi_{0}^{\rm out} (t,{\bf y}_{1/2\pm}) \cdot {\bf n}_0.
\end{equation}
\end{subequations}
Moreover, in the following, we extend the derivative condition for the mole fraction 
\begin{equation}\label{ES3-10}
\lim_{z \rightarrow \pm \infty} \left( \partial_z \chi_{1}^{v,\rm in} (t,z,{\bf s}) ~ {\bf n}_0 + \nabla_{\Gamma} \chi_{0}^{v,\rm in} (t,z,{\bf s}) \right) = \nabla \chi_{0}^{v,\rm out} (t,{\bf y}_{1/2\pm}),
\end{equation}
and similarly for the velocity
\begin{subequations}\label{ES3-11}
\begin{equation}\label{ES3-11a}
\lim_{z \rightarrow \pm \infty} \left( \partial_{z} {\rm {\bf v}}_{1}^{\rm in} (t,z,{\bf s}) \otimes {\bf n}_{0} + \nabla_{\Gamma} {\rm {\bf v}}_{0}^{\rm in} (t,z,{\bf s}) \right) = \nabla {\rm {\bf v}}_{0}^{\rm out} (t,{\bf y}_{1/2\pm}),
\end{equation}
\begin{equation}\label{ES3-11b}
\lim_{z \rightarrow \pm \infty} \left( \partial_{z} {\rm {\bf v}}_{1}^{\rm in} (t,z,{\bf s}) \cdot {\bf n}_{0} + \nabla_{\Gamma} \cdot {\rm {\bf v}}_{0}^{\rm in} (t,z,{\bf s}) \right) = \nabla \cdot {\rm {\bf v}}_{0}^{\rm out} (t,{\bf y}_{1/2\pm}).
\end{equation}
\end{subequations}

\subsection{Outer Expansions}\label{Sec:4-1}
We substitute the outer expansions for the unknowns into the equations \eqref{ES2-1}, \eqref{ES2-3}, \eqref{ES2-4}, \eqref{ES2-6} and \eqref{ES2-8}. Then collecting terms for the leading order leads to
\begin{equation}\label{ES3-12}
P'(\phi_{0}^{out}) = 0.
\end{equation}
\begin{equation}\label{ES3-13}
\partial_{t} \rho \left( \phi_{0}^{\rm out} \right) + \nabla \cdot \left( \rho \left( \phi_{0}^{\rm out} \right) {\rm {\bf v}}_{0}^{\rm out} \right) = 0,
\end{equation}
\begin{equation}\label{ES3-14}
\partial_{t} \left( \rho (\phi_{0}^{\rm out}) \chi_{0}^{v,\rm out} \right) + \nabla \cdot \left( \rho (\phi_{0}^{\rm out}) \chi_{0}^{v,\rm out} {\rm {\bf v}}_{0}^{\rm out} \right) = \nabla \cdot \left( {\rm D}^{v} \rho (\phi_{0}^{\rm out}) \nabla \chi_{0}^{v,\rm out} \right),
\end{equation}
\begin{equation}\label{ES3-15}
\partial_{t} \left( \rho (\phi_{0}^{\rm out}) {\rm {\bf v}}_{0}^{\rm out} \right) + \nabla \cdot \left( \rho (\phi_{0}^{\rm out}) {\rm {\bf v}}_{0}^{\rm out} \otimes {\rm {\bf v}}_{0}^{\rm out} \right) = - \nabla p_{0}^{\rm out} + \nabla \cdot \boldsymbol{\tau}_{0}^{\rm out} + \rho (\phi_{0}^{\rm out}) {\rm {\bf g}},
\end{equation}
\begin{equation}\label{ES3-16}
\begin{aligned}
\partial_{t} \left( \rho(\phi^{out}_{0}) u_{0}^{out} \right) + \nabla \cdot \left( \rho(\phi^{out}_{0}) h_{0}^{out} {\rm {\bf v}}_{0}^{out} \right) = \nabla \cdot \left( k(\phi_{0}^{out}) \nabla {\rm T}_{0}^{out} \right) & + {\rm {\bf v}}^{out}_{0} \cdot \nabla p^{out}_{0} \\
 & + \boldsymbol{\tau}^{out}_{0} : \nabla {\rm {\bf v}}^{out}_{0}.
\end{aligned}
\end{equation}
Note that the equation \eqref{ES3-12} has three solutions as: $\phi_{0}^{out} = 0, ~ 1/2, ~ \text{or} ~ 1$. Following similar arguments as in Ref. \cite{vanNoorden2011}, the solution $\phi_{0}^{out} = 1/2$ is unstable as $P''(1/2) < 0$, whereas the first and the last solution are stable since $P''(0)$ and $P''(1)$ are positive. In view of this, we see that in the limit $\zeta \rightarrow 0$ one obtains the solutions $\phi_{0}^{out} = 0$ and $\phi_{0}^{out} = 1$. Now we associate the notations $\Omega_{0}^{l}(t)$ and $\Omega_{0}^{g}(t)$ as the (time dependent) sub-domains of $\Omega$ where $\phi_{0}^{out}$ is 1 and 0, respectively. Then, the governing equations for $\phi_{0}^{out} = 1$ in $\Omega_{0}^{l}(t)$ are
\begin{equation}\label{ES3-17}
\partial_{t} \rho_{l,0}^{out} + \nabla \cdot \left( \rho_{l,0}^{out} {\rm {\bf v}}_{l,0}^{out} \right) = 0,
\end{equation}
\begin{equation}\label{ES3-18}
\partial_{t} \left( \rho_{l,0}^{out} \chi_{l,0}^{v,out} \right) + \nabla \cdot \left( \rho_{l,0}^{out} \chi_{l,0}^{v,out} {\rm {\bf v}}_{l,0}^{out} \right) = \nabla \cdot \left( {\rm D}_{l}^{v} \rho_{l,0}^{out} \nabla \chi_{l,0}^{v,out} \right),
\end{equation}
\begin{equation}\label{ES3-19}
\begin{aligned}
    \partial_{t} \left( \rho_{l,0}^{out} {\rm {\bf v}}_{l,0}^{out} \right) + \nabla \cdot \left( \rho_{l,0}^{out} {\rm {\bf v}}_{l,0}^{out} \otimes {\rm {\bf v}}_{l,0}^{out} \right) = & - \nabla p_{l,0}^{out} + \nabla \cdot \left\{ \mu_{l} \left( \nabla {\rm {\bf v}}_{l,0}^{out} + \nabla {\rm {\bf v}}_{l,0}^{out^{T}} \right) \right\} \\
     & + \nabla \cdot \left( \xi_{l} \nabla \cdot {\rm {\bf v}}_{l,0}^{out} {\rm {\bf I}} \right) + \rho_{l,0}^{out} {\rm {\bf g}},
\end{aligned}
\end{equation}
\begin{equation}\label{ES3-20}
\begin{aligned}
\partial_{t} \left( \rho_{l,0}^{out} u_{l,0}^{out} \right) + \nabla \cdot \left( \rho_{l,0}^{out} h_{l,0}^{out} {\rm {\bf v}}_{l,0}^{out} \right) & = \nabla \cdot \left( {\rm k}_l \nabla {\rm T}_{l,0}^{out} \right) + {\rm {\bf v}}^{out}_{l,0} \cdot \nabla p^{out}_{l,0} \\
& + \left\{ \mu_{l} \left( \nabla {\rm {\bf v}}_{l,0}^{out} + \nabla {\rm {\bf v}}_{l,0}^{out^{T}} \right) + \xi_{l} \nabla \cdot {\rm {\bf v}}_{l,0}^{out} {\rm {\bf I}} \right\} : \nabla {\rm {\bf v}}^{out}_{l,0}.
\end{aligned}
\end{equation}
Next the governing equations for $\phi_{0}^{out} = 0$ in $\Omega_{0}^{g}(t)$ are
\begin{equation}\label{ES3-21}
\partial_{t} \rho_{g,0}^{out} + \nabla \cdot \left(  \rho_{g,0}^{out} {\rm {\bf v}}_{g,0}^{out} \right) = 0,
\end{equation}
\begin{equation}\label{ES3-22}
\partial_{t} \left( \rho_{g,0}^{out} \chi_{g,0}^{v,out} \right) + \nabla \cdot \left( \rho_{g,0}^{out} \chi_{g,0}^{v,out} {\rm {\bf v}}_{g,0}^{out} \right) = \nabla \cdot \left( {\rm D}_{g}^{v} \rho_{g,0}^{out} \nabla \chi_{g,0}^{v,out} \right),
\end{equation}
\begin{equation}\label{ES3-23}
\begin{aligned}
    \partial_{t} \left( \rho_{g,0}^{out} {\rm {\bf v}}_{g,0}^{out} \right) + \nabla \cdot \left( \rho_{g,0}^{out} {\rm {\bf v}}_{g,0}^{out} \otimes {\rm {\bf v}}_{g,0}^{out} \right) = & - \nabla p_{g,0}^{out} + \nabla \cdot \left\{ \mu_{g} \left( \nabla {\rm {\bf v}}_{g,0}^{out} + \nabla {\rm {\bf v}}_{g,0}^{out^{T}} \right) \right\} \\
     &  + \nabla \cdot \left(\xi_{g} \nabla \cdot {\rm {\bf v}}_{g,0}^{out} {\rm {\bf I}} \right) + \rho_{g,0}^{out} {\rm {\bf g}},
\end{aligned}
\end{equation}
\begin{equation}\label{ES3-24}
\begin{aligned}
\partial_{t} \left( \rho_{g,0}^{out} u_{g,0}^{out} \right) + \nabla \cdot \left( \rho_{g,0}^{out} h_{g,0}^{out} {\rm {\bf v}}_{g,0}^{out} \right) & = \nabla \cdot \left( {\rm k}_g \nabla {\rm T}_{g,0}^{out} \right) + {\rm {\bf v}}^{out}_{g,0} \cdot \nabla p^{out}_{g,0} \\
& + \left\{ \mu_{g} \left( \nabla {\rm {\bf v}}_{g,0}^{out} + \nabla {\rm {\bf v}}_{g,0}^{out^{T}} \right) + \xi_{g} \nabla \cdot {\rm {\bf v}}_{g,0}^{out} {\rm {\bf I}} \right\} : \nabla {\rm {\bf v}}^{out}_{g,0}.
\end{aligned}
\end{equation}
Here, we note that the above equations for the liquid and gas regions are exactly the same model equations as described in the sharp-interface formulation (see \hyperref[Sec:2]{Section 2}).

\subsection{Inner Expansions}\label{Sec:4-2}
We now apply the inner expansions and the matching conditions to the phase-field model to recover the interface conditions on the evolving interface.

\subsubsection{Phase-field Equation}\label{Sec:4-2-1}
The dominating $\mathcal{O}(\zeta^{-2})$ terms in the phase-field equation \eqref{ES2-1} are
\begin{equation}\label{ES3-25}
L^2 \partial^{2}_{z} \phi_{0}^{in} - P'(\phi_{0}^{in}) = 0.
\end{equation}
Using the matching condition (\ref{ES3-9a}), we have 
\[
\lim_{z \rightarrow -\infty} \phi_{0}^{in}(t,z,{\bf s}) = 1 \qquad \text{and} \qquad \lim_{z \rightarrow + \infty} \phi_{0}^{in}(t,z,{\bf s}) = 0.
\]
Moreover, when $\zeta \rightarrow 0$, the moving interface is defined by $\phi_{0}^{in}(t,0,{\bf s}) = 1/2$. Hence, multiplying equation
\eqref{ES3-25} by $\partial_{z} \phi_{0}^{in}$, integrating the result with respect to $z$ and using the matching conditions fulfilled by $\phi_{0}^{in}$ (see (\ref{ES3-9})) and the specific form of $P(\phi)$, we have
\begin{equation}\label{ES3-26}
\partial_{z} \phi_{0}^{in} = - \frac{\sqrt{2}}{L} \phi_{0}^{in} \left( 1 - \phi_{0}^{in} \right).
\end{equation}
Using $\phi_{0}^{in}(t,0,{\bf s}) = 1/2$, the solution becomes
\begin{equation}\label{ES3-27}
\phi_{0}^{in}(t,z,{\bf s}) = \phi_{0}^{in}(z) = \frac{1}{1 + e^{\frac{\sqrt{2}z}{L}}} = \frac{1}{2} \left( 1 - \tanh \left( \frac{z}{\sqrt{2}L} \right) \right).
\end{equation}
Now for the $\mathcal{O}(\zeta^{-1})$ terms, we have
\begin{equation}\label{ES3-28}
\begin{aligned}
\left( P''(\phi_{0}^{in}) - L^2 \partial^{2}_{z} \right) \phi_{1}^{in} = & \left( L^2 \kappa_{0} + \frac{L^2}{\gamma/\nu} \rho (\phi_{0}^{in}) v_{n,0} \right) \partial_z \phi_{0}^{in} - \frac{L^2}{\gamma/\nu} \rho (\phi_{0}^{in}) \partial_z \left( \phi_{0}^{in} {\rm {\bf v}}_{0}^{in} \right) \cdot {\rm {\bf n}}_{0} \\
& - \frac{\sqrt{2} L}{\gamma/\nu} \phi_{0}^{in} (1 - \phi_{0}^{in}) f(\chi_{0}^{v,in}, {\rm T}_{0}^{in}).
\end{aligned}
\end{equation}
Let us assume that $\partial_z {\rm {\bf v}}_{0}^{in} \cdot {\bf n}_{0} = 0$, see \hyperref[Re:1]{Remark 4.1} for more details. In view of this, equation \eqref{ES3-28} reduces to
\begin{dmath}\label{ES3-29}
\left( P''(\phi_{0}^{in}) - L^2 \partial^{2}_{z} \right) \phi_{1}^{in} = \left\{ L^2 \kappa_{0} + \frac{L^2}{\gamma/\nu} \rho (\phi_{0}^{in}) \left( v_{n,0} - {\rm {\bf v}}_{0}^{in} \cdot {\rm {\bf n}}_{0} \right) \right\} \partial_z \phi_{0}^{in} - \frac{\sqrt{2} L}{\gamma/\nu} \phi_{0}^{in} (1 - \phi_{0}^{in}) f(\chi_{0}^{v,in}, {\rm T}_{0}^{in}).
\end{dmath}
We consider that $\mathcal{L} (\phi_{0}^{in}) = P''(\phi_{0}^{in}) - L^2 \partial^{2}_{z}$ is an operator, which is applied to $\phi_{1}^{in}$ on the left-hand side of the above equation. Note that $\mathcal{L}$ is a Fredholm operator of index zero, hence the equation (\ref{ES3-29}) possesses a unique solution if and only if the right-hand side, denoted by $\mathcal{F}(\phi_{0}^{in})$, is orthogonal to the kernel of $\mathcal{L}$. From (\ref{ES3-25}), one can see that the kernel of $\mathcal{L}$ contains $\partial_z \phi_{0}^{in}$. Now, since $v_{n,0}$, $\kappa_0$, $\chi_{0}^{v,in}$, and ${\rm T}_{0}^{in}$ are independent of $z$ (these will be shown in the forthcoming subsection), the solvability condition implies
\begin{equation}\label{ES3-30}
\int_{-\infty}^{+\infty} \mathcal{F} \left( \phi_{0}^{in} \right) \partial_{z} \phi_{0}^{in} dz = 0,
\end{equation}
where $\mathcal{F} \left( \phi_{0}^{in} \right) = \left\{ L^2 \kappa_{0} + \frac{L^2}{\gamma/\nu} \rho (\phi_{0}^{in}) \left( v_{n,0} - {\rm {\bf v}}_{0}^{in} \cdot {\rm {\bf n}}_{0} \right) \right\} \partial_z \phi_{0}^{in} - \frac{\sqrt{2} L}{\gamma/\nu} \phi_{0}^{in} (1 - \phi_{0}^{in}) f(\chi_{0}^{v,in}, {\rm T}_{0}^{in}).$ We also assume that the densities $\rho_l$ and $\rho_g$ are independent of $z$. Then, equation \eqref{ES3-30} can be simplified as
\begin{equation}\label{ES3-31}
\frac{\left( \rho_{l,0}^{in} + \rho_{g,0}^{in} \right)}{2} \left( v_{n,0} - {\rm {\bf v}}_{0}^{in} \cdot {\rm {\bf n}}_{0} \right) = - \left( \frac{\gamma}{\nu} \kappa_{0} + f(\chi_{0}^{v,in}, {\rm T}_{0}^{in}) \right).
\end{equation}
Now upon applying the matching conditions (\ref{ES3-9}), the above relation (\ref{ES3-31}) reduces to
\begin{equation}\label{ES3-32}
\frac{\left( \rho_{l,0}^{out} + \rho_{g,0}^{out} \right)}{2} \left( v_{n,0} - {\rm {\bf v}}_{l,0}^{out} \cdot {\rm {\bf n}}_{0} \right) = - \left( \frac{\gamma}{\nu} \kappa_{0} + f(\chi_{l,0}^{v,out}, {\rm T}_{l,0}^{out}) \right),
\end{equation}
which is the same as the reaction rate (\ref{ES1-26}) on the sharp interface $\Gamma_{lg} (t)$ except for the curvature term. Although it is a common property of both two-phase flow processes and in particular evaporation processes to be affected by curvature, the term $\frac{\gamma}{\nu}\kappa_0$ arises due to the Allen-Cahn formulation and will not necessarily represent a physically correct motion. Such curvature-driven motion can however be minimized by choosing the fraction $\frac{\gamma}{\nu}$ small.

\begin{remark}\label{Re:1}
Note that the present contribution deals with a compressible phase-field model and we are mainly interested in a diffusion-dominated process, where viscous terms play an important role. However, Dreyer \emph{et al.} \cite{Dreyer2014} showed that it is not admissible to have phase change for compressible phase-field models when including viscous terms. Hence, to still include phase change in the presence of viscous terms, we have chosen  $\partial_z {\rm {\bf v}}_{0}^{in} \cdot {\bf n}_{0} = 0$ in the previous derivation, which corresponds to relaxing the compressibility of the fluids. This assumption is not necessary when considering the mass conservation.
\end{remark}

\subsubsection{Mass Conservation Equations}\label{Sec:4-2-2}
The dominating $\mathcal{O}(\zeta^{-1})$ terms arising from inserting the inner expansions into equation \eqref{ES2-3} are
\begin{equation}\label{ES3-33}
v_{n,0} \partial_z \rho \left( \phi_{0}^{in} \right) - \partial_z \left( \rho \left( \phi_{0}^{in} \right) {\rm {\bf v}}_{0}^{in} \right) \cdot {\rm {\bf n}}_{0} = 0.
\end{equation}
Integrating with respect to $z$ and using matching conditions (\ref{ES3-9}) leads to
\begin{equation}\label{ES3-34}
v_{n,0} \left( \rho_{g,0}^{out} - \rho_{l,0}^{out} \right) = \left( \rho_{g,0}^{out} {\rm {\bf v}}_{g,0}^{out} - \rho_{l,0}^{out} {\rm {\bf v}}_{l,0}^{out} \right) \cdot {\rm {\bf n}}_{0},
\end{equation}
which is equivalent to the transmission condition (\ref{ES1-2}) on the sharp interface $\Gamma_{lg} (t)$. Next, the dominating $\mathcal{O}(\zeta^{-2})$ term arising from inserting the inner expansions into equation \eqref{ES2-4} is
\begin{equation}\label{ES3-35}
\partial_z \left( {\rm D}^{v} \rho(\phi_{0}^{in}) \partial_z \chi_{0}^{v,in} \right) = 0.
\end{equation}
Now integrating \eqref{ES3-35} with respect to $z$, we get
\begin{equation}\label{ES3-36}
{\rm D}^{v} \rho(\phi_{0}^{in}) \partial_z \chi_{0}^{v,in} = C_{1} = constant.
\end{equation}
By determining the constant through matching conditions (\ref{ES3-9}), the above equation gives that $\partial_z\chi_{0}^{v,\text{in}} = 0$. Again from equation \eqref{ES2-4}, the $\mathcal{O}(\zeta^{-1})$ terms give
\begin{equation}\label{ES3-37}
\begin{aligned}
v_{n,0} \partial_z \left( \rho(\phi_{0}^{in}) \chi_{0}^{v,in} \right) - \partial_z \left( \rho(\phi_{0}^{in}) \chi_{0}^{v,in} {\rm {\bf v}}_{0}^{in} \right) \cdot {\rm {\bf n}}_{0} = - \partial_z \left( {\rm D}^{v} \phi_{1}^{in} \rho'(\phi_{0}^{in}) \partial_z \chi_{0}^{v,in} \right)\\
- \partial_z \left( {\rm D}^{v} \rho(\phi_{0}^{in}) \partial_z \chi_{1}^{v,in} \right) - \partial_z \left( {\rm D}^{v} \rho(\phi_{0}^{in}) \nabla_{\Gamma} \chi_{0}^{v,in} \right) \cdot {\rm {\bf n}}_{0} - \nabla_{\Gamma} \cdot \left( {\rm D}^{v} \rho(\phi_{0}^{in}) \partial_z \chi_{0}^{v,in} {\rm {\bf n}}_{0} \right).
\end{aligned}
\end{equation}
Using the fact that $\partial_z\chi_{0}^{v,\text{in}} = 0$ in \eqref{ES3-37}, and then integrating \eqref{ES3-37} along with the matching conditions (\ref{ES3-9}) and (\ref{ES3-10}), we obtain
\begin{equation}\label{ES3-38}
\begin{aligned}
    v_{n,0} \left( \rho_{g,0}^{out} \chi_{g,0}^{v,out} - \rho_{l,0}^{out} \chi_{l,0}^{v,out} \right) & = \left( \rho_{g,0}^{out} \chi_{g,0}^{v,out} {\rm {\bf v}}_{g,0}^{out} - \rho_{l,0}^{out} \chi_{l,0}^{v,out} {\rm {\bf v}}_{l,0}^{out} \right) \cdot {\rm {\bf n}}_{0} \\
    & - \left( {\rm D}_{g}^{v} \rho_{g,0}^{out} \nabla \chi_{g,0}^{v,out} - {\rm D}_{l}^{v} \rho_{l,0}^{out} \nabla \chi_{l,0}^{v,out}\right) \cdot {\rm {\bf n}}_{0}.
\end{aligned}
\end{equation}
Note that the above condition is the same as the one given in (\ref{ES1-5}) on the sharp interface $\Gamma_{lg} (t)$.

\subsubsection{Momentum Balance Equation}\label{Sec:4-2-3}
The dominating $\mathcal{O}(\zeta^{-2})$ terms in the momentum equation \eqref{ES2-6} are
\begin{equation}\label{ES3-39}
\begin{aligned}
L \sigma \partial_{z} \left( \partial_{z} \phi_{0}^{in} \right)^2 {\bf n}_{0} - \partial_{z} \left( \mu \left( \phi_{0}^{in} \right) \left( \partial_{z} {\rm {\bf v}}_{0}^{in} \otimes {\bf n}_{0} + \left( \partial_{z} {\rm {\bf v}}_{0}^{in} \otimes {\bf n}_{0} \right)^{T} \right) {\bf n}_{0} \right) & \\
- \partial_{z} \left( \xi \left( \phi_{0}^{in} \right) \partial_{z} {\rm {\bf v}}_{0}^{in} \cdot {\bf n}_{0} {\rm {\bf I}} \right) {\bf n}_{0} & = 0.
\end{aligned}
\end{equation}
Simple algebra leads to
\begin{equation}\label{ES3-40}
\begin{aligned}
L \sigma \partial_{z} \left( \partial_{z} \phi_{0}^{in} \right)^2 {\bf n}_{0} - \partial_{z} \left( \mu (\phi_{0}^{in}) \partial_{z} {\rm {\bf v}}_{0}^{in} \right) - \partial_{z} \left( \mu (\phi_{0}^{in}) \left( \partial_{z} {\rm {\bf v}}_{0}^{in} \otimes {\bf n}_{0} \right)^{T} \right) {\bf n}_{0} & \\
- \partial_{z} \left( \xi (\phi_{0}^{in}) \partial_{z} {\rm {\bf v}}_{0}^{in} \cdot {\bf n}_{0} \right) {\bf n}_{0} & = 0.
\end{aligned}
\end{equation}
We introduce the tangent vector ${\bf t}_{0}$. Note that ${\bf t}_{0}$ is independent of $z$, and that ${\bf n}_{0} \cdot {\bf t}_{0} = 0$. Taking the dot product of equation \eqref{ES3-40} with ${\bf t}_{0}$ leads to
\begin{equation}\label{ES3-41}
\partial_{z} \left( \mu(\phi_{0}^{in}) \partial_{z} {\rm {\bf v}}_{0}^{in} \cdot {\bf t}_{0} \right) = 0.
\end{equation}
Integrating with respect to $z$ and using matching conditions (\ref{ES3-9}) results in
\begin{equation}\label{ES3-42}
\mu(\phi_{0}^{in}) \partial_{z} {\rm {\bf v}}_{0}^{in} \cdot {\bf t}_{0} = constant.
\end{equation}
Since $\lim_{z \rightarrow \pm \infty} \partial_{z} {\rm {\bf v}}_{0}^{in} = 0$, the integration constant must be 0. Now as $\mu(\phi_{0}^{in}) \neq 0$, equation \eqref{ES3-42} implies
\begin{equation}\label{ES3-43}
\partial_{z} {\rm {\bf v}}_{0}^{in} \cdot {\bf t}_{0} = 0.
\end{equation}
Integrating again with respect to $z$ and using matching condition (\ref{ES3-9a}) gives
\begin{equation}\label{ES3-44}
{\rm {\bf v}}_{l,0}^{out} \left( t, {\rm {\bf y}}_{1/2} \right) \cdot {\bf t}_{0} = {\rm {\bf v}}_{g,0}^{out} \left( t, {\rm {\bf y}}_{1/2} \right) \cdot {\bf t}_{0},
\end{equation}
which is the same as the condition (\ref{ES1-8}) at the sharp interface $\Gamma_{lg} (t)$. We assume that $\partial_{z} {\rm {\bf v}}_{0}^{in} \cdot {\bf n}_{0} = 0$ and now also derive that $\partial_{z} {\rm {\bf v}}_{0}^{in} \cdot {\bf t}_{0} = 0$. Hence, from these two we get $\partial_{z} {\rm {\bf v}}_{0}^{in} = 0$. Next, for the $\mathcal{O}(\zeta^{-1})$ terms, we have
\begin{equation}\label{ES3-45}
\begin{aligned}
- v_{n,0} \partial_{z} \left( \rho \left( \phi_{0}^{in} \right) {\rm {\bf v}}_{0}^{in} \right) & + \partial_{z} \left( \rho \left( \phi_{0}^{in} \right) {\rm {\bf v}}_{0}^{in} \otimes {\rm {\bf v}}_{0}^{in} \right) {\bf n}_{0} = - \partial_{z} p_{0}^{in} {\bf n}_{0} \\
 & + \partial_{z} \left( \mu \left( \phi_{0}^{in} \right) \left( \partial_{z} {\rm {\bf v}}_{1}^{in} \otimes {\bf n}_{0} + \left( \partial_{z} {\rm {\bf v}}_{1}^{in} \otimes {\bf n}_{0} \right)^{T} \right) {\bf n}_{0} \right) \\
 & + \partial_{z} \left( \mu \left( \phi_{0}^{in} \right) \left( \nabla_{\Gamma} {\rm {\bf v}}_{0}^{in} + \left( \nabla_{\Gamma} {\rm {\bf v}}_{0}^{in} \right)^{T} \right) {\bf n}_{0} \right) \\
 & + \nabla_{\Gamma} \cdot \left( \mu \left( \phi_{0}^{in} \right) \left( \partial_{z} {\rm {\bf v}}_{0}^{in} \otimes {\bf n}_{0} + \left( \partial_{z} {\rm {\bf v}}_{0}^{in} \otimes {\bf n}_{0} \right)^{T} \right) \right) \\
 & + \partial_{z} \left( \mu' \left( \phi_{0}^{in} \right) \phi_{1}^{in} \left( \partial_{z} {\rm {\bf v}}_{0}^{in} \otimes {\bf n}_{0} + \left( \partial_{z} {\rm {\bf v}}_{0}^{in} \otimes {\bf n}_{0} \right)^{T} \right) {\bf n}_{0} \right) \\
 & + \partial_{z} \left( \xi \left( \phi_{0}^{in} \right) \partial_{z} {\rm {\bf v}}_{1}^{in} \cdot {\bf n}_{0} {\rm {\bf I}} \right) {\bf n}_{0} + \partial_{z} \left( \xi \left( \phi_{0}^{in} \right) \nabla_{\Gamma} \cdot {\rm {\bf v}}_{0}^{in} {\rm {\bf I}} \right) {\bf n}_{0} \\
 & + \nabla_{\Gamma} \cdot \left( \xi \left( \phi_{0}^{in} \right) \partial_{z} {\rm {\bf v}}_{0}^{in} \cdot {\bf n}_{0} {\rm {\bf I}} \right) + \partial_{z} \left( \xi' \left( \phi_{0}^{in} \right) \phi_{1}^{in} \partial_{z} {\rm {\bf v}}_{0}^{in} \cdot {\bf n}_{0} {\rm {\bf I}} \right) {\bf n}_{0} \\
 & + L \sigma \left( \partial_{z} \phi_{0}^{in} \right)^{2} \kappa_{0} {\bf n}_{0} - L \sigma \partial_{z} \left( \partial_{z} \phi_{1}^{in} \partial_{z} \phi_{0}^{in} \right) {\bf n}_{0}.
\end{aligned}
\end{equation}
Now, since $\partial_z {\rm {\bf v}}_{0}^{in} = 0$, then equation \eqref{ES3-45} reduces to
\begin{equation}\label{ES3-46}
\begin{aligned}
- v_{n,0} \partial_{z} \left( \rho \left( \phi_{0}^{in} \right) {\rm {\bf v}}_{0}^{in} \right) & + \partial_{z} \left( \rho \left( \phi_{0}^{in} \right) {\rm {\bf v}}_{0}^{in} \otimes {\rm {\bf v}}_{0}^{in} \right) {\bf n}_{0} = - \partial_{z} p_{0}^{in} {\bf n}_{0} \\
 & + \partial_{z} \left( \mu \left( \phi_{0}^{in} \right) \left( \partial_{z} {\rm {\bf v}}_{1}^{in} \otimes {\bf n}_{0} + \left( \partial_{z} {\rm {\bf v}}_{1}^{in} \otimes {\bf n}_{0} \right)^{T} \right) {\bf n}_{0} \right) \\
 & + \partial_{z} \left( \mu \left( \phi_{0}^{in} \right) \left( \nabla_{\Gamma} {\rm {\bf v}}_{0}^{in} + \left( \nabla_{\Gamma} {\rm {\bf v}}_{0}^{in} \right)^{T} \right) {\bf n}_{0} \right) \\
 & + \partial_{z} \left( \xi \left( \phi_{0}^{in} \right) \partial_{z} {\rm {\bf v}}_{1}^{in} \cdot {\bf n}_{0} {\rm {\bf I}} \right) {\bf n}_{0} + \partial_{z} \left( \xi \left( \phi_{0}^{in} \right) \nabla_{\Gamma} \cdot {\rm {\bf v}}_{0}^{in} {\rm {\bf I}} \right) {\bf n}_{0} \\
 & + L \sigma \left( \partial_{z} \phi_{0}^{in} \right)^{2} \kappa_{0} {\bf n}_{0} - L \sigma \partial_{z} \left( \partial_{z} \phi_{1}^{in} \partial_{z} \phi_{0}^{in} \right) {\bf n}_{0}.
\end{aligned}
\end{equation}
Integrating \eqref{ES3-46} with respect to $z$ and upon applying matching conditions (\ref{ES3-9}) and \eqref{ES3-11}, leads to
\begin{dmath}\label{ES3-47}
- v_{n,0} \left[ \rho \left( \phi_{0}^{out} \right) {\rm {\bf v}}_{0}^{out} \right]_{-\infty}^{+\infty} + \left[ \rho \left( \phi_{0}^{out} \right) {\rm {\bf v}}_{0}^{out} \otimes {\rm {\bf v}}_{0}^{out} \right]_{-\infty}^{+\infty} {\bf n}_{0} - L \sigma \left( \int_{-\infty}^{+\infty} \left( \partial_{z} \phi_{0}^{in} \right)^{2} dz \right) \kappa_{0} {\bf n}_{0} = - \left[ p_{0}^{out} \right]_{-\infty}^{+\infty} {\bf n}_{0} + \left[ \mu \left( \phi_{0}^{out} \right) \left( \nabla {\rm {\bf v}}_{0}^{out} + \left( \nabla {\rm {\bf v}}_{0}^{out} \right)^{T} \right) + \xi \left( \phi_{0}^{out} \right) \nabla \cdot {\rm {\bf v}}_{0}^{out} {\rm {\bf I}} \right]_{-\infty}^{+\infty}  {\bf n}_{0}.
\end{dmath}
Note that the value of the integral $\int_{-\infty}^{+\infty} \left( \partial_{z} \phi_{0}^{in} \right)^{2} dz$ is $\frac{\sqrt{2}}{6 L}$. Then, upon substituting the limit values for $z \rightarrow \pm \infty$ in the above relation leads to the condition (\ref{ES1-7}) on the sharp interface $\Gamma_{lg} (t)$.

\subsubsection{Energy Balance Equation}\label{Sec:4-2-4}
The dominating $\mathcal{O}(\zeta^{-2})$ terms in the energy equation \eqref{ES2-8} are
\begin{equation}\label{ES3-48}
\begin{aligned}
\partial_{z} \left( {\rm k}(\phi_{0}^{in}) \partial_{z} {\rm T}_{0}^{in} \right) & + \mu(\phi_{0}^{in}) \left( \partial_z {\rm {\bf v}}^{in}_{0} \otimes {\bf n}_{0} + \left( \partial_z {\rm {\bf v}}^{in}_{0} \otimes {\bf n}_{0} \right)^{T} \right) : \left( \partial_z {\rm {\bf v}}^{in}_{0} \otimes {\bf n}_{0} \right) \\
& + \xi(\phi_{0}^{in}) \left( \partial_z {\rm {\bf v}}^{in}_{0} \cdot {\bf n}_{0} \right) {\bf I} : \left( \partial_z {\rm {\bf v}}^{in}_{0} \otimes {\bf n}_{0} \right) = 0.
\end{aligned}
\end{equation}
Now upon using the fact that $\partial_z {\rm {\bf v}}_{0}^{in} = 0$ in the above equation \eqref{ES3-48}, leads to
\begin{equation}\label{ES3-49}
\partial_{z} \left( {\rm k}(\phi_{0}^{in}) \partial_{z} {\rm T}_{0}^{in} \right) = 0.
\end{equation}
Integrating with respect to $z$ and using matching conditions (\ref{ES3-9b}) and the fact that ${\rm k}(\phi_{0}^{in}) > 0$, we obtain
\begin{equation}\label{ES3-50}
\partial_{z} {\rm T}_{0}^{in} = 0.
\end{equation}
Integrating again with respect to $z$ and using matching conditions (\ref{ES3-9a}), we get
\begin{equation}\label{ES3-51}
{\rm T}_{g,0}^{out} \left( t, {\rm {\bf y}}_{1/2} \right) = {\rm T}_{l,0}^{out} \left( t, {\rm {\bf y}}_{1/2} \right),
\end{equation}
which is the same as the condition (\ref{ES1-16}) on the sharp interface $\Gamma_{lg} (t)$. Taking advantage of $\partial_{z} {\rm T}_{0}^{in} = 0$, the $\mathcal{O}(\zeta^{-1})$ terms satisfy
\begin{equation}\label{ES3-52}
\begin{aligned}
- v_{n,0} \partial_{z} \left( \rho(\phi_{0}^{in}) u_{0}^{in} \right) + \partial_{z} \left( \rho(\phi_{0}^{in}) h_{0}^{in} {\rm {\bf v}}_{0}^{in} \right) \cdot {\bf n}_{0} = \partial_{z} \left( {\rm k}(\phi_{0}^{in}) \partial_{z} {\rm T}_{1}^{in} \right) + {\rm {\bf v}}^{in}_{0} \cdot \partial_z p^{in}_{0} {\bf n}_{0} \\
+ \left\{ \mu(\phi_{0}^{in}) \left( \partial_z {\rm {\bf v}}^{in}_{0} \otimes {\bf n}_{0} + \left( \partial_z {\rm {\bf v}}^{in}_{0} \otimes {\bf n}_{0} \right)^{T} \right) + \xi(\phi_{0}^{in}) \left( \partial_z {\rm {\bf v}}^{in}_{0} \cdot {\bf n}_{0} \right) {\bf I} \right\} : \left( \partial_z {\rm {\bf v}}^{in}_{1} \otimes {\bf n}_{0} \right) \\
+ \left\{ \mu(\phi_{0}^{in}) \left( \partial_z {\rm {\bf v}}^{in}_{1} \otimes {\bf n}_{0} + \left( \partial_z {\rm {\bf v}}^{in}_{1} \otimes {\bf n}_{0} \right)^{T} \right) + \xi(\phi_{0}^{in}) \left( \partial_z {\rm {\bf v}}^{in}_{1} \cdot {\bf n}_{0} \right) {\bf I} \right\} : \left( \partial_z {\rm {\bf v}}^{in}_{0} \otimes {\bf n}_{0} \right) \\
+ \left\{ \mu(\phi_{0}^{in}) \left( \partial_z {\rm {\bf v}}^{in}_{0} \otimes {\bf n}_{0} + \left( \partial_z {\rm {\bf v}}^{in}_{0} \otimes {\bf n}_{0} \right)^{T} \right) + \xi(\phi_{0}^{in}) \left( \partial_z {\rm {\bf v}}^{in}_{0} \cdot {\bf n}_{0} \right) {\bf I} \right\} : \nabla_{\Gamma} {\rm {\bf v}}^{in}_{0} \\
+ \left\{ \mu(\phi_{0}^{in}) \left( \nabla_{\Gamma} {\rm {\bf v}}^{in}_{0} + \left( \nabla_{\Gamma} {\rm {\bf v}}^{in}_{0} \right)^{T} \right) + \xi(\phi_{0}^{in}) \left( \nabla_{\Gamma} \cdot {\rm {\bf v}}^{in}_{0} \right) {\bf I} \right\} : \left( \partial_z {\rm {\bf v}}^{in}_{0} \otimes {\bf n}_{0} \right).
\end{aligned}
\end{equation}
Again using $\partial_z {\rm {\bf v}}_{0}^{in} = 0$ in the previous equation \eqref{ES3-52} leads to
\begin{equation}\label{ES3-53}
- v_{n,0} \partial_{z} \left( \rho(\phi_{0}^{in}) u_{0}^{in} \right) + \partial_{z} \left( \rho(\phi_{0}^{in}) h_{0}^{in} {\rm {\bf v}}_{0}^{in} \right) \cdot {\bf n}_{0} = \partial_{z} \left( {\rm k}(\phi_{0}^{in}) \partial_{z} {\rm T}_{1}^{in} \right) + {\rm {\bf v}}^{in}_{0} \cdot  {\bf n}_{0} ~ \partial_z p^{in}_{0}.
\end{equation}
Integrating with respect to $z$ and using the matching conditions (\ref{ES3-9}), we get
\begin{equation}\label{ES3-54}
\begin{aligned}
- v_{n,0} \left( \rho_{g,0}^{out} u_{g,0}^{out} - \rho_{l,0}^{out} u_{l,0}^{out} \right) + \left( \rho_{g,0}^{out} h_{g,0}^{out} {\rm {\bf v}}_{g,0}^{out} - \rho_{l,0}^{out} h_{l,0}^{out} {\rm {\bf v}}_{l,0}^{out} \right) \cdot {\bf n}_{0} \\
= \left( {\rm k}_{g} \nabla {\rm T}_{g,0}^{out} - {\rm k}_{l} \nabla {\rm T}_{l,0}^{out} \right) \cdot {\bf n}_{0} + \left( p_{g,0}^{out} {\rm {\bf v}}_{g,0}^{out} - p_{l,0}^{out} {\rm {\bf v}}_{l,0}^{out} \right) \cdot {\bf n}_{0}.
\end{aligned}
\end{equation}
Using $\dot{m} = \rho_{\alpha} \left( {\rm {\bf v}}_{\alpha} \cdot {\bf n} - v_{n} \right)$, the above condition can also be written as
\begin{equation}\label{ES3-55}
\left( {\rm k}_{g} \nabla {\rm T}_{g,0}^{out} - {\rm k}_{l} \nabla {\rm T}_{l,0}^{out} \right) \cdot {\bf n}_{0} = \dot{m}_{0} \left( u_{g,0}^{out} - u_{l,0}^{out} \right).
\end{equation}
Note that this is the flux condition (\ref{ES1-18}) on the sharp interface $\Gamma_{lg} (t)$. Similar rewriting as discussed in \hyperref[Sec:2-3]{Section 2.3} is of course also applicable here.

\section{Numerical Behavior}\label{Sec:5}
We consider two numerical examples to analyze the applicability and the potential of the phase-field formulation developed in this study. In both cases we assume a liquid droplet sitting on a solid wall. The focus of the present study is to analyze the evaporation from a porous medium, consequently the flow velocity is smaller, hence we neglect the inertial term in the momentum balance equation \eqref{ES2-6} as well as the body force term for simplicity. Next, we also neglect the last three terms from the energy balance equation \eqref{ES2-8}, since these are having less significant influence than the other terms. One can realize that by analyzing the order of the corresponding terms considering the process of interest. In the first example, we consider that the domain is closed for fluid flow, however, the second example allows flow from the left to the right of the domain. All model equations are implemented in DuMux \cite{Koch2021}, utilizing Newton’s method as the nonlinear solver and the implicit Euler method for temporal discretization. Spatial discretization is achieved using the finite volume method on a regular rectangular grid, employing Dune-SPGrid \cite{Nolte2011}. A cell-centered approach is used for the phase field, vapor concentration, pressure and temperature unknowns, with control volumes matching the grid cells and degrees of freedom positioned at their centers. Fluxes across control-volume faces are approximated using the two-point flux approximation (TPFA) approach, \emph{i.e.,} using the values from the two adjacent cells. For the velocities, a staggered discretization is employed, with separate degrees of freedom for each velocity component placed at the edges between grid cells, and control volumes centered around them, for more details refer \cite{Koch2021}. As we are mainly interested in the evaporation process and consequently on the phase change phenomena, hence, we use the IAPWS Formulation 1995 \cite{IAPWS1995} for the fluid phase properties. Detailed description of the domain properties along with the initial and boundary conditions are described in the following subsections. 

\subsection{Example I: No External Flow Field}\label{Sec:5-1}
\begin{figure}[h!]
    \centering
    \begin{overpic}[scale = 0.095]{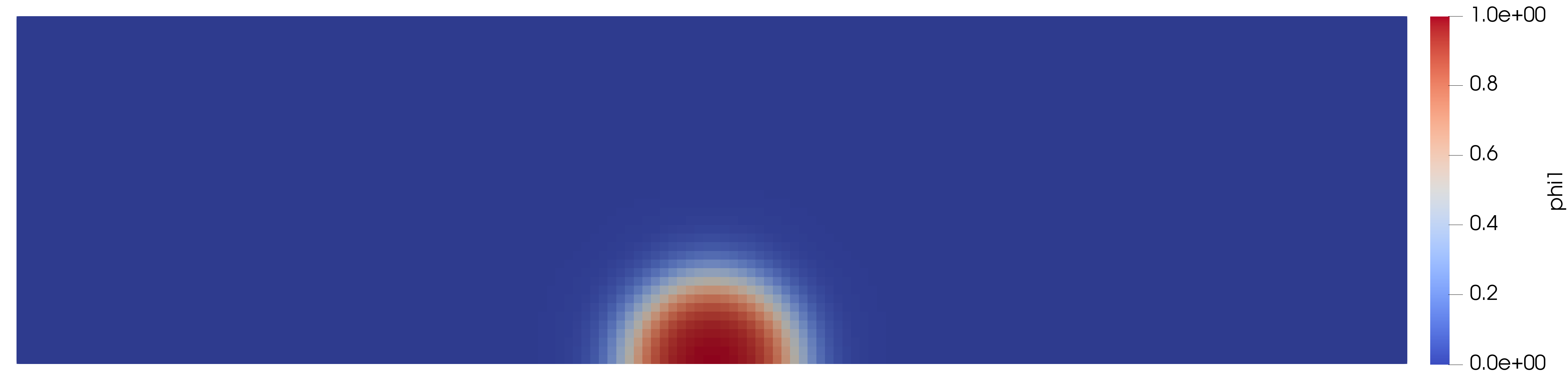}
    \put(-13,14){\bf No flow}
    \put(-13,10){\bf No flux}
    \put(101,14){\bf No flow}
    \put(101,10){\bf No flux}
    \put(32,-3){\bf No flow, No flux}
    \put(32,24){\bf No flow, No flux}
    \end{overpic}
    \caption{Initial phase distribution for the no external flow field along with the boundary conditions.}
    \label{noflow_geometry}
\end{figure}
\begin{figure}[h!]
    \centering
    \includegraphics[scale=0.1]{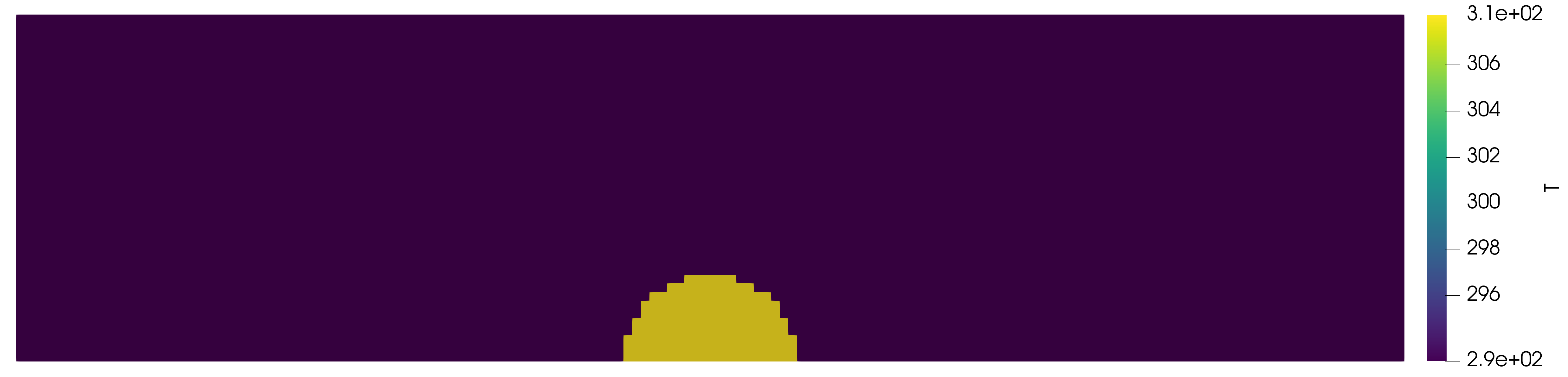}
    \caption{Initial temperature distribution for both the cases of induced and no external flow field.}
    \label{temperature_initial}
\end{figure}

\noindent We consider a two-dimensional rectangular channel of length $2~mm$ and width $0.5~mm$. We initialize the fluid distribution as shown in \hyperref[noflow_geometry]{Figure 4}, where $\phi = 1$ (red) represents liquid phase and $\phi = 0$ (blue) corresponds to gas phase with a diffusive layer between them. The liquid droplet is of diameter $0.25~mm$. Further, we assume that initially the liquid droplet is having higher temperature ($308.15~K$) relative to the surrounding gas phase temperature ($293.15~K$), see \hyperref[temperature_initial]{Figure 5}. All the boundaries of the channel are impermeable, hence we incorporate no flow condition at all the boundaries. In addition, we also have no-flux boundary conditions for density, vapor concentration, temperature and phase field at all the boundaries. Although the system is closed, we choose the temperature and vapor concentration in such a way that the evaporation rate, as defined in \hyperref[Sec:2-4]{Section 2.4}, remains positive, thereby initiating droplet evaporation.

\begin{figure}[h!]
    \centering
    \includegraphics[scale=0.1]{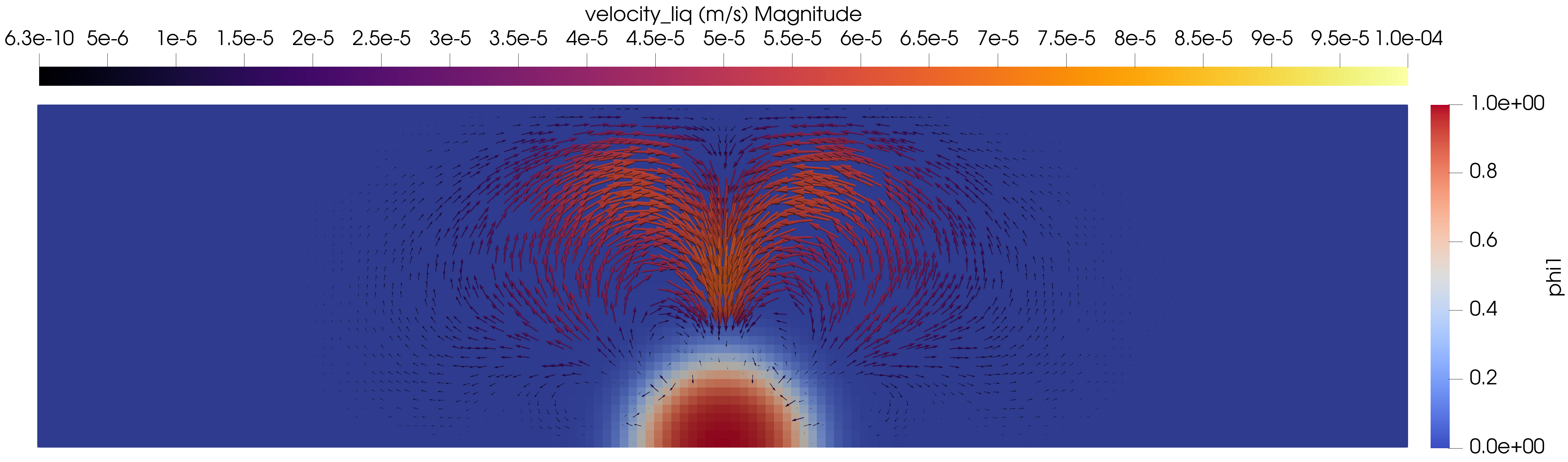}
    \caption{Velocity field overlaid on phase distribution at time $t = 40~ns$ for the no external flow field.}
    \label{noflow_velocity}
\end{figure}
\begin{figure}[h!]
    \centering
    \includegraphics[scale=0.1]{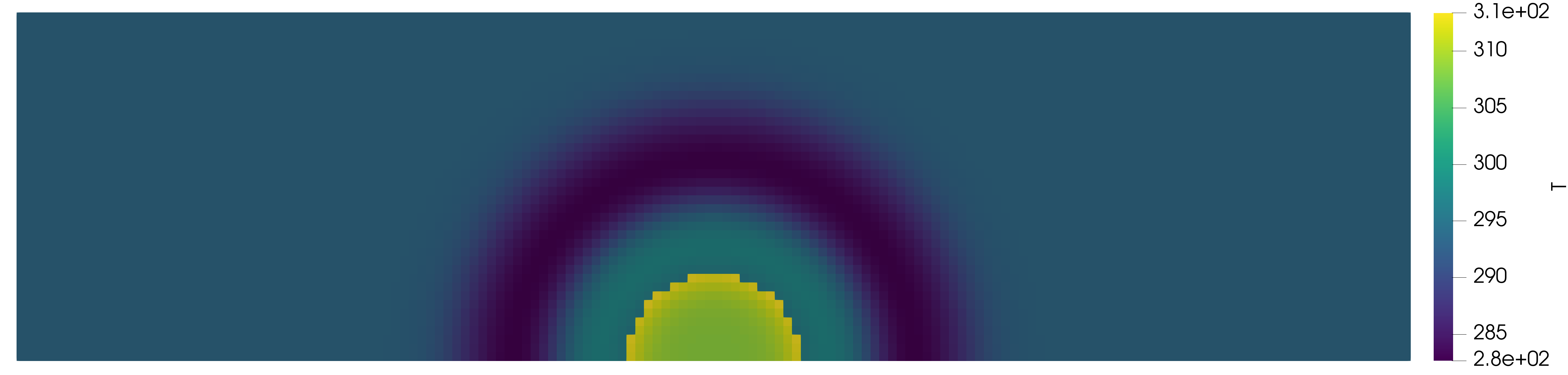}
    \caption{Temperature distribution around the liquid droplet at time $t = 40~ns$ for the no external flow field.}
    \label{noflow_temperature}
\end{figure}
\hyperref[noflow_velocity]{Figure 6} shows the phase distribution along with the overlaid velocity field at time $t = 40~ns$. One can notices a relatively higher velocity magnitude around the droplet, specifically at the top of the droplet. Although the domain has no externally imposed flow field, one can still observe flow around the droplet, see \hyperref[noflow_velocity]{Figure 6}. However, the flow magnitude is quite small, which eventually suggests a diffusion dominated process. 
Initially starting with a hot drop, one can clearly distinguish the evaporation related cooling effect near the droplet at time $t = 40~ns$, see \hyperref[noflow_temperature]{Figure 7}. Further, looking into the temperature distribution, one can certainly observe that there is a steep temperature gradient around the droplet. Since we initialize the temperature with a sharp transition on the droplet interface, a diffusive layer in the temperature distribution between the drop and the cooling region appears, which is an artifact of the phase-field approach. However this is less significant than the advantage this formulation possesses to model moving interfaces. Further, in \hyperref[noflow_flow_intphi]{Figure 11}, we show the changes in the size of the liquid droplet against time. 

\subsection{Example II: Induced External Flow Field}\label{Sec:5-2}
\begin{figure}[h!]
    \centering
    \begin{overpic}[scale = 0.1]{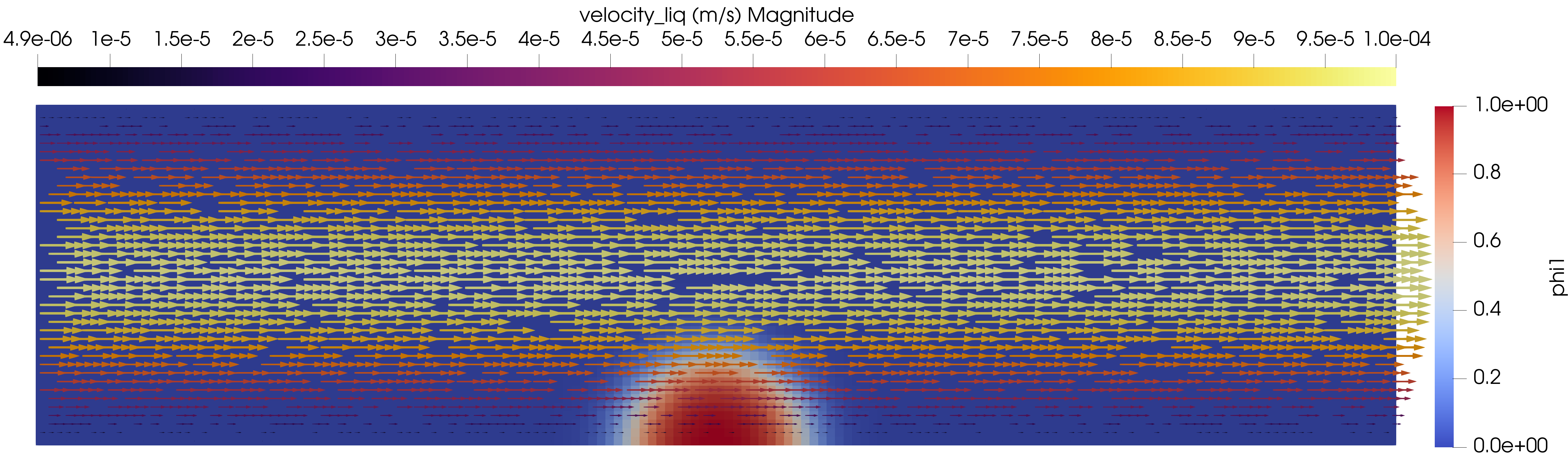}
    \put(-10,12){\bf Inflow}
    \put(102,14){\bf Fixed}
    \put(102,10){\bf Pressure}
    \put(32,-2){\bf No flow, No flux}
    \put(32,30){\bf No flow, No flux}
    \end{overpic}
    \caption{Initial phase distribution for the induced external flow field along with the boundary conditions.}
    \label{flow_geometry}
\end{figure}
\noindent For this example, again we consider a $2~mm \times 0.5~mm$ rectangular channel. We initialize the fluid distribution as shown in \hyperref[flow_geometry]{Figure 8}, where $\phi = 1$ (red) represents liquid phase and $\phi = 0$ (blue) corresponds to gas phase with a diffusive layer between them. The liquid droplet is of diameter $0.25~mm$. Further, the initial velocity field (fully developed and parabolic) overlaid on phase distribution can be observed in \hyperref[flow_geometry]{Figure 8}. In a similar way like the previous example, we again consider a hot liquid droplet as shown in \hyperref[temperature_initial]{Figure 5}. The top and bottom boundaries of the channel are impermeable, and no-flux boundary conditions for density, vapor concentration, temperature and phase field are employed at these locations. On the left boundary, we impose an inlet velocity and apply Dirichlet conditions for all other variables: the vapor concentration and phase-field are set to zero, and the temperature is $298.15 K$. On the other hand, at the right boundary, we maintain a fixed pressure of $10^5 Pa$ and apply homogeneous Neumann (no flux) conditions for all other variables. 

\begin{figure}[h!]
    \centering
    \includegraphics[scale=0.1]{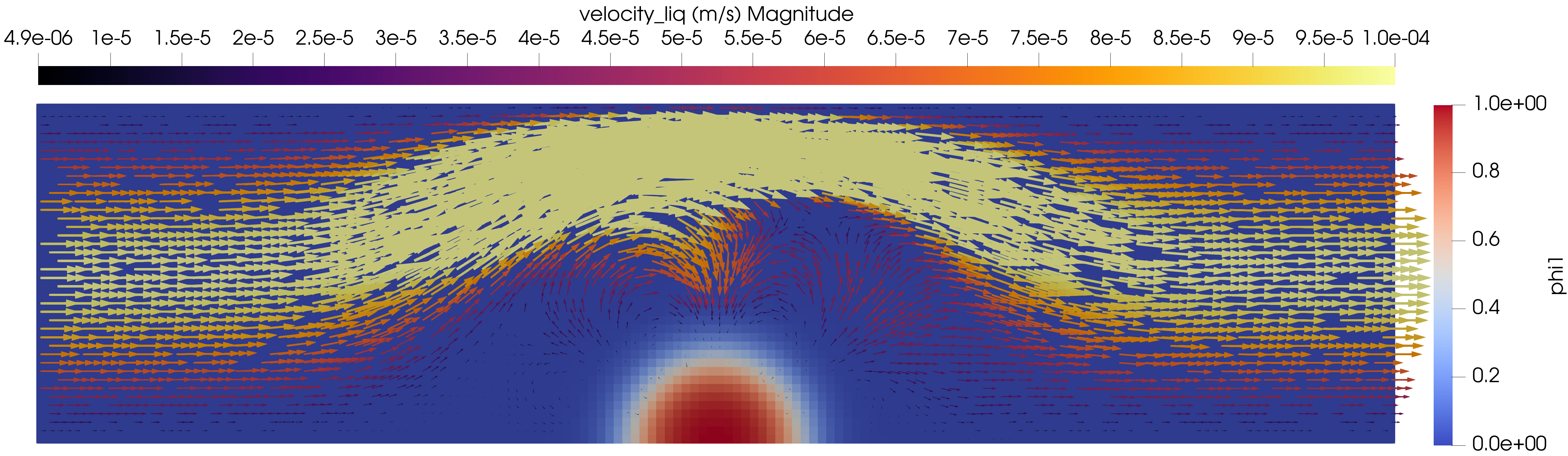}
    \caption{Velocity field overlaid on phase distribution at time $t = 26~ns$ for the induced external flow field.}
    \label{flow_velocity}
\end{figure}
\begin{figure}[h!]
    \centering
    \includegraphics[scale=0.1]{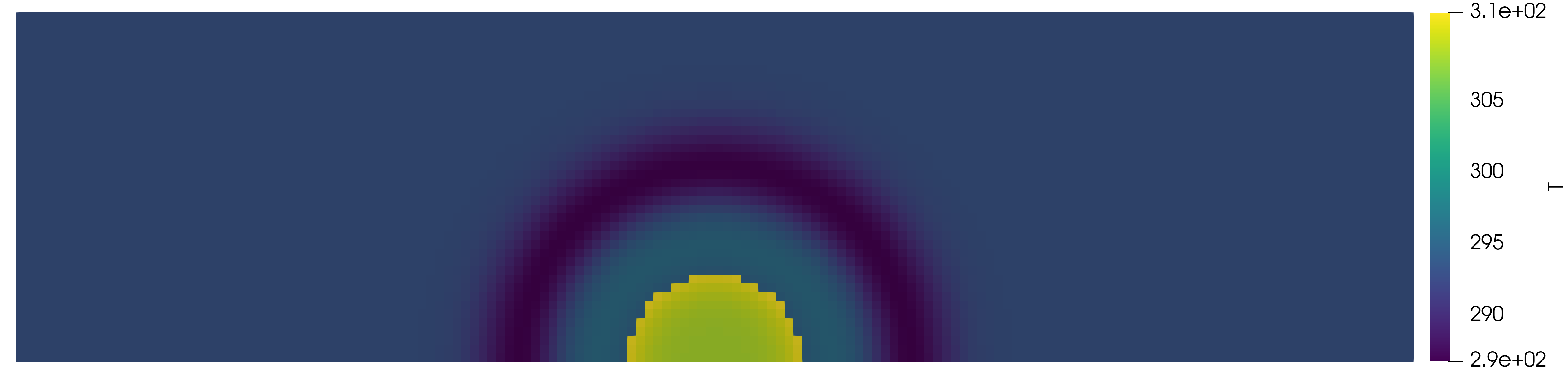}
    \caption{Temperature distribution around the liquid droplet at time $t = 26~ns$ for the induced external flow field.}
    \label{flow_temperature}
\end{figure}
\hyperref[flow_velocity]{Figure 9} shows the phase distribution along with the overlaid velocity field at time $t = 26~ns$. A notable change in the flow field surrounding the droplet is readily apparent. Looking into the flow field one can interpret that the droplet behaves as a slight obstacle to the flow situation. Due to the difference in density and viscosity of the droplet relative to the neighbouring gas phase, we notice a relatively small flow magnitude around the droplet. Although very near to the droplet (specifically in the diffuse interface layer) we notice diffusion phenomena, but away from the droplet the induced external flow field dominates the process. Again considering a hot drop initially, we observe evaporation related cooling effect around the droplet at time $t = 26~ns$ (see \hyperref[flow_temperature]{Figure 10}), which we also observed in the no external flow case (\hyperref[noflow_temperature]{Figure 7}). Furthermore, we note a similar kind of temperature gradient around the droplet as in the case of no external flow field. 

\begin{figure}[h!]
\centering
\includegraphics[scale=0.25]{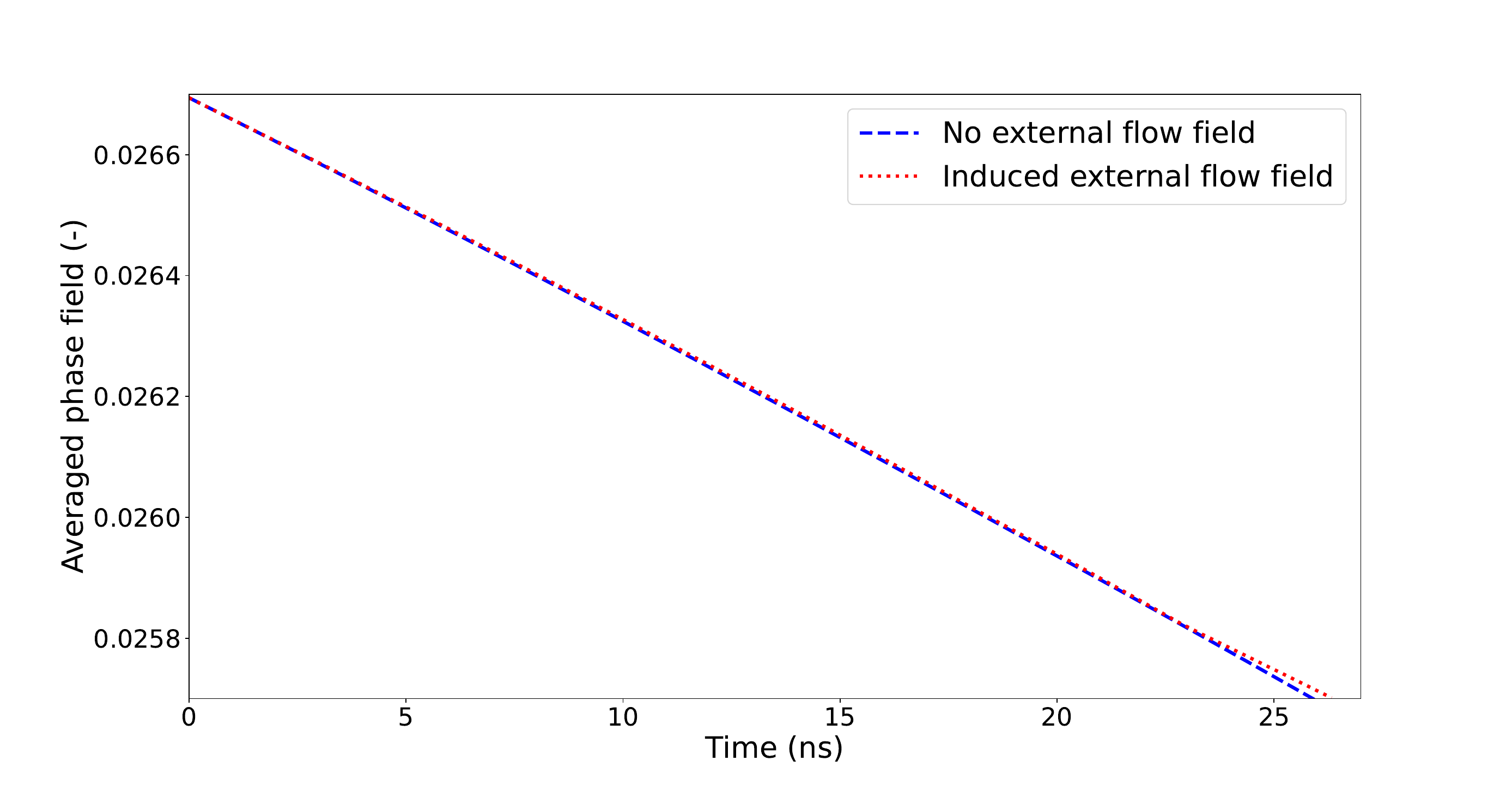}
\caption{Changes in the averaged phase field against time for no external and induced external flow field.}
\label{noflow_flow_intphi}
\end{figure} 

We compute the averaged phase field by integrating the phase field over the domain and then dividing by the volume, resulting in a quantity that varies only with time. The changes in the averaged phase field over time are depicted in \hyperref[noflow_flow_intphi]{Figure 11}. We can certainly observe that the averaged phase field decreases with time, indicating that the liquid droplet is shrinking, which implies its evaporation into vapor. After $26~ns$, approximately $3.7\%$ of the droplet has evaporated under both no external flow and induced external flow conditions (see \hyperref[noflow_flow_intphi]{Figure 11}). However, there are minor deviations in the evaporation behavior between the two conditions around $25~ns$. The simulation's wall time is restricted due to the challenges the Newton method faces in dealing with the complex coupled system analyzed in this study.


\section{Upscaling: Periodic Homogenization}\label{Sec:6}
The phase-field model described in \hyperref[Sec:3]{Section 3} approximates the evolution of two phases in an arbitrary domain with an arbitrary number of phase transitions. Moreover, the phase-field model is a pore-scale model, which is defined in a very complex domain. Now in such a scenario, one is mainly interested in the averaged behavior of the system. As a consequence, we apply two-scale asymptotic expansion techniques to find effective equations and parameters at a larger scale, which does not experience the micro-scale properties but still includes their effect.

\begin{figure}[h!]
    \centering
    \begin{overpic}[width=0.85\textwidth]{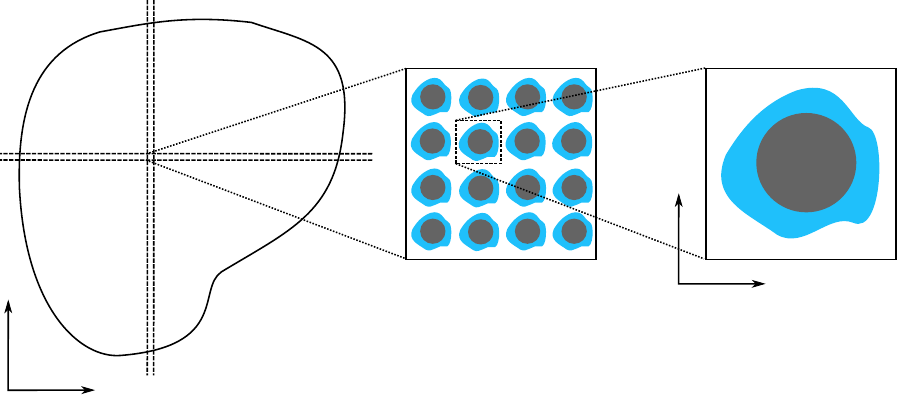}
    \put(11.5,0){\color{black} $x_1$}
    \put(-1,12){\color{black} $x_2$}
    \put(86,12){\color{black} $y_1$}
    \put(74,24){\color{black} $y_2$}
    \end{overpic}
    \caption{Sketch of a porous medium, depicting liquid as the wetting phase. Void space filled with gas and liquid are marked with white and blue, respectively; non-reactive solid phase (grains) is in grey.}
    \label{F4}
\end{figure}
Now we consider an idealized configuration where solid (grains) and pores are periodically distributed, which is usually referred to as a periodic porous medium. We assume that the phase-field formulation given in equations \eqref{ES2-1}, \eqref{ES2-3}, \eqref{ES2-4}, \eqref{ES2-6}, and \eqref{ES2-8} is defined in such a medium. Taking into account the effective behavior of the system, in the following we derive a macro (Darcy) scale model. Specifically, we consider a domain $\mathbb{D}$, which contains periodically distributed solid phase (grains), see Figure \ref{F4}. Note that $\mathbb{D}$ is the union of the grain space and the void (pore) space for a porous medium. We assume that the grains are impermeable and does not react with the fluid phases; however, grains can exchange heat with the fluids to balance the total energy of the system. Therefore, the phase-field formulation described in equations \eqref{ES2-1}, \eqref{ES2-3}, \eqref{ES2-4}, \eqref{ES2-6}, and \eqref{ES2-8} is only defined in the pore space of the domain $\mathbb{D}$. Moreover, we assume that the pore space in $\mathbb{D}$ is connected.

\subsection{Two Scales}\label{Sec:6-1}
The parameter $\ell$ stands for a characteristic length scale at the pore (micro) scale, \emph{e.g.,} the width of the rightmost box in \hyperref[F4]{Figure 12}, and $L$ is a reference length scale at the macro level, \emph{e.g.,} the width of the domain $\mathbb{D}$ or of the Darcy scale, as typically considered for homogenization \cite{Davit2013,Hornung1996}. Then we define $\varepsilon = \ell/L$, which represents the ratio between the micro and macro scales and consequently gives us the scale separation. Generally, $\ell$ is much smaller than $L$, thus $\varepsilon$ is a small number. Precisely, micro-scale variable ${\bf y} \in (0, \ell)^d$ and macro-scale variable ${\bf x} \in (0, L)^d$, where $d$ refers to the spatial dimension, which is either 2 or 3. Since non-dimensionalization is an essential step for upscaling, we start with non-dimensionalization of the domain. We define the non-dimensional space variables as follows: $\hat{\rm {\bf x}} = {\rm {\bf x}}/L$ and $\hat{\rm {\bf y}} = {\rm {\bf y}}/\ell$. Hereafter, we only consider the dimensionless domain and for notational convenience we drop the hat.

\begin{figure}[h!]
    \centering
    \begin{overpic}[width=0.4\textwidth]{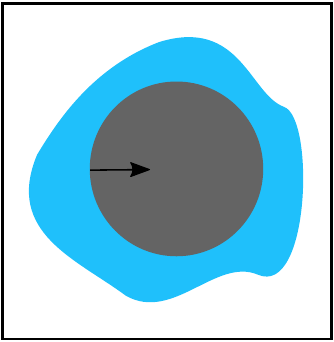}
    \put(52,48){\color{black} $S$}
    \put(29,53){\color{black} ${\bf n}_P$}
    \put(23,65){\color{black} $W$}
    \put(10,70){\color{black} $G$}
    \put(40,25){\color{black} $\Gamma_P$}
    \put(45,2){\color{black} $\partial Y$}
    \end{overpic}
    \caption{Sketch of a local pore $Y = [0,1]^d$. Gas region (white) is $G$, liquid region (blue) is $W$, and solid/grain portion (gray) is $S$. The normal vector at the internal boundary $\Gamma_P$ is ${\bf n}_P$. The outer boundary of the local pore, $\partial Y$, is marked in black.}
    \label{F5}
\end{figure}
In the following subsection we rewrite the governing equations in non-dimensional form. In connection to this we first define a local unit cell $Y = [0,1]^{d}$ (see \hyperref[F5]{Figure 13}). Moreover, we assume that the local variable ${\bf y} \in (0,1)^{d}$ represent points within $Y$. Note that the local cell contains the gas part $G$, liquid part $W$, and the solid (grain) part $S$, as depicted in \hyperref[F5]{Figure 13}. Therefore, locally the phase-field model is defined in the void (pore) space $P = G \cup W$, and $S$ defines the solid. The boundary $\Gamma_P$ denotes the (stationary) interface between the domain of the phase-field model and the solid. While $\partial Y$ defines the outer boundary of the unit cell $Y$, where we will use the periodic boundary conditions to decouple the unit cells from each other. In the following whenever we mention internal boundaries, we generally mean $\Gamma_P$.

Now to differentiate between the two scales in the considered model, we define ${\bf x}$ as the macro scale variable, and ${\bf y}$ as the micro scale variable. These are connected through the relation ${\bf y} = \varepsilon^{-1} {\bf x}$. In other words, ${\bf x}$ only sees the macro-scale characteristics, whereas the zooomed-in ${\bf y}$ observes the micro-scale oscillations in a single cell. Therefore, for each macro-scale point ${\bf x}$, one can determine a unit cell, having its own local variable ${\bf y}$. Note that the domain of the phase-field model is the union of all the local pore space $P$, scaled by $\varepsilon$, while the domain for the solid heat transport is the union of all local solids $S$. Hence, the domains depend on $\varepsilon$ and can be expressed as
\[ \Omega^{\varepsilon} = \cup_{w \in \mathbb{W}_{\mathbb{D}}} \left\{ \varepsilon (w + P) \right\} \quad \text{and} \quad \Omega_{S}^{\varepsilon} = \cup_{w \in \mathbb{W}_{\mathbb{D}}} \left\{ \varepsilon (w + S) \right\}, \]
where $\mathbb{W}_{\mathbb{D}}$ is a subset of $\mathbb{Z}^{d}$ satisfying $\mathbb{D} = \cup_{w \in \mathbb{W}_{\mathbb{D}}} \left\{ \varepsilon (w + Y) \right\}$, which is the entire domain, see the leftmost box in \hyperref[F4]{Figure 12}. To show the dependence on $\varepsilon$, we use $\varepsilon$ as a superscript. The union of all internal solid-pore space boundaries $\Gamma_P$ is expressed as
\[ \Gamma^{\varepsilon} = \cup_{w \in \mathbb{W}_{\mathbb{D}}} \left\{ \varepsilon (w + \Gamma_P) \right\}.  \]

\subsection{Non-Dimensional Model Equations}\label{Sec:6-2}
We first non-dimensionalize the phase-field formulation described by the equations \eqref{ES2-1}--\eqref{ES2-10} and then the heat conduction equation (\ref{ES1-23}) for the solid phase, to identify the dominating terms and consequently the significant ones for the upscaling. To ensure that we are in the range of Darcy's law, \emph{i.e.,} the momentum equation \eqref{ES2-6} becomes a Darcy-like equation, we assume the typical velocity, pressure, viscosity, density etc. accordingly. Moreover, we also ensure that the diffuse interface (\emph{i.e.}, the transition between liquid and gas) stays within a local pore. Non-dimensional variables and quantities are denoted with a hat and are defined as follows
\[
\begin{array}{ l@{\hskip 0.5in}  l@{\hskip 0.5in}  l@{\hskip 0.5in} l }
\hat{t} = t/t_{\rm ref}, & \hat{\rho} = \rho/\rho_{\rm ref}, & \hat{\rm {\bf v}} = {\rm {\bf v}}/{\rm V}_{\rm ref}, & \hat{p} = p/p_{\rm ref}, \\
\hat{\mu} = \mu/\mu_{\rm ref}, & \hat{\xi} = \xi/\xi_{\rm ref}, & \hat{\rm D}^{v} = {\rm D}^{v}/{\rm D}_{\rm ref}, & \hat{\lambda} = \lambda/\lambda_{\rm ref}, \\
\hat{\sigma} = \sigma/\sigma_{\rm ref}, & \hat{\rm {\bf g}} = {\rm {\bf g}}/{\rm g}_{\rm ref}, & \hat{\gamma} = \gamma/\gamma_{\rm ref}, & \hat{\nu} = \nu/\nu_{\rm ref}, \\
\hat{u} = u/u_{\rm ref}, & \hat{h} = h/u_{\rm ref}, & \hat{\rm k} = {\rm k}/{\rm k}_{\rm ref}, & \hat{\rm T} = {\rm T}/{\rm T}_{\rm ref}, \\
\hat{\rho}_{S} = \rho_{S}/\rho_{\rm ref}, & \hat{\rm T}_{S} = {\rm T}_{S}/{\rm T}_{\rm ref}, & \hat{\rm k}_{S} = {\rm k}_{S}/{\rm k}_{\rm ref}, & \hat{\rm c}_{p,S} = {\rm c}_{p,S}/{\rm c}_{\rm ref}.
\end{array}
\]
Note that the reference quantities are related to each other through various non-dimensional numbers. One can consider the size of these non-dimensional numbers depending on the regime of interest. Here in the present case, we are interested in Darcy's regime where the diffusion time scale dominates the advection and reaction time scales. It is well-known that the Darcy's law is valid when fluid flow is laminar and the pressure drop dominates the flow behavior, which corresponds to the Reynolds and Euler numbers being
\[
{\rm Re} = \frac{\rho_{\rm ref} {\rm V}_{\rm ref} L}{\mu_{\rm ref}} = \mathcal{O}(\varepsilon^0) \quad \text{and} \quad {\rm Eu} = \frac{p_{\rm ref}}{\rho_{\rm ref} {\rm V}^{2}_{\rm ref}} = \mathcal{O}(\varepsilon^{-2}),
\]
respectively. Advection, diffusion and reaction time scales are identified as 
\[
t_{\rm adv} = \frac{L}{{\rm V}_{\rm ref}}, \quad t_{\rm diff} = \frac{L^{2}}{{\rm D}_{\rm ref}} \quad \text{and} \quad t_{\rm react} = \frac{\rho_{\rm ref} L}{R_{\rm ref}},
\] 
respectively. It is assumed that diffusion is more important than advection and reaction, which corresponds to the following scaling of the P\'eclet and Damköhler numbers
\[
{\rm Pe} = \frac{t_{\rm diff}}{t_{\rm adv}} = \frac{L {\rm V}_{\rm ref}}{{\rm D}_{\rm ref}} = \mathcal{O}(\varepsilon) \quad \text{and} \quad {\rm Da} = \frac{t_{\rm adv}}{t_{\rm react}} = \frac{R_{\rm ref}}{\rho_{\rm ref} {\rm V}_{\rm ref}} = \mathcal{O}(\varepsilon^{0}).
\]
We assume that the observation time scale $t_{\rm ref}$ is equal to the diffusion time scale $t_{\rm diff}$. The width of the diffuse interface layer is proportional to, but smaller than, the pore size $\ell$, which corresponds to the Cahn number being
\[
{\rm Ch} = \frac{\lambda_{\rm ref}}{L} = \mathcal{O}(\varepsilon).
\]
Therefore, $\hat{\lambda}$ is independent of $\varepsilon$ and a small number. Consequently, the width of the diffuse interface layer is small compared to the pore size, but nevertheless a micro-scale quantity which will not approach its sharp-interface limit when $\varepsilon$ approaches zero.

Now, to capture the influence of heat conduction on the leading order equations at the Darcy scale, we choose the Prandtl number ${\rm Pr} = u_{\rm ref} \mu_{\rm ref} / {\rm k}_{\rm ref} {\rm T}_{\rm ref} = \mathcal{O}(\varepsilon^{0})$ and Schmidt number ${\rm Sc} = \mu_{\rm ref} / \rho_{\rm ref} {\rm D}_{\rm ref} = \mathcal{O}(\varepsilon^{0})$. In order to keep the gravity effect on the momentum equation, we choose the Froude number ${\rm Fr} = {\rm V}_{\rm ref} / \sqrt{{\rm g}_{\rm ref} L} = \mathcal{O}(\varepsilon)$. Note that, in the present case, surface tension forces are dominating the viscous drag forces, \emph{i.e.,} Capillary number ${\rm Ca} = \mu_{\rm ref} {\rm V}_{\rm ref} / \sigma_{\rm ref} = \mathcal{O}(\varepsilon^{2})$. Next, we assume that the dynamic and the kinematic viscosity are of the same order, \emph{i.e.,} $\xi_{\rm ref} / \mu_{\rm ref} = \mathcal{O}(\varepsilon^{0})$. Finally, we choose $\gamma_{\rm ref} / \nu_{\rm ref} \mu_{\rm ref} = \mathcal{O}(\varepsilon)$ in order to have a consistent scaling of the terms in the phase-field equation. For readability, we assume that the non-dimensional numbers that are $\mathcal{O}(\varepsilon^{n})$ to be exactly equal to $\varepsilon^{n}$.

In the following we will only use non-dimensional variables, so we skip the hat on all variables. Consequently, the dimensionless model reads
\begin{equation}\label{ES4-1}
\begin{aligned}
\rho^{\varepsilon} \left\{ \partial_t \phi^{\varepsilon} + \varepsilon \nabla \cdot \left( \phi^{\varepsilon} {\rm {\bf v}}^{\varepsilon}  \right) \right\} & = \varepsilon^{2} \frac{\gamma}{\nu} \left\{ \nabla^{2} \phi^{\varepsilon} - \varepsilon^{-2} \lambda^{-2} P'(\phi^{\varepsilon}) \right\} \\
& - \sqrt{2} \lambda^{-1} \phi^{\varepsilon} (1 - \phi^{\varepsilon}) f(\chi^{v \varepsilon}, {\rm T}^{\varepsilon}) \qquad \text{in } \Omega^{\varepsilon},
\end{aligned}
\end{equation}
\begin{equation}\label{ES4-2}
\partial_t \rho^{\varepsilon} + \varepsilon \nabla \cdot \left( \rho^{\varepsilon} {\rm {\bf v}}^{\varepsilon} \right) = 0 \qquad \text{in } \Omega^{\varepsilon},
\end{equation}
\begin{equation}\label{ES4-3}
\partial_t \left( \rho^{\varepsilon} \chi^{v \varepsilon} \right) + \varepsilon \nabla \cdot \left( \rho^{\varepsilon} \chi^{v\varepsilon} {\rm {\bf v}}^{\varepsilon} \right) = \nabla \cdot \left( {\rm D}^{v} \rho^{\varepsilon} \nabla \chi^{v \varepsilon} \right) \qquad \text{in } \Omega^{\varepsilon},
\end{equation}
\begin{equation}\label{ES4-4}
\begin{aligned}
\partial_t \left( \rho^{\varepsilon} {\rm {\bf v}}^{\varepsilon} \right) + \varepsilon \nabla \cdot \left( \rho^{\varepsilon} {\rm {\bf v}}^{\varepsilon} \otimes {\rm {\bf v}}^{\varepsilon} \right) & = - \varepsilon^{-1} \nabla p^{\varepsilon} + \varepsilon \nabla \cdot \left\{ \mu^{\varepsilon} \left( \nabla {\rm {\bf v}}^{\varepsilon} + \nabla {\rm {\bf v}}^{{\varepsilon}^{\rm T}} \right) + \xi^{\varepsilon} \left( \nabla \cdot {\rm {\bf v}}^{\varepsilon} \right) {\bf I} \right\} \\
 & - \varepsilon^{2} \nabla \cdot \left( \lambda \sigma \nabla \phi^{\varepsilon} \otimes \phi^{\varepsilon} \right) + \varepsilon^{-1} \rho^{\varepsilon} {\rm {\bf g}} \qquad \text{in } \Omega^{\varepsilon},
\end{aligned}
\end{equation}
\begin{equation}\label{ES4-5}
\begin{aligned}
\partial_t \left( \rho^{\varepsilon} u^{\varepsilon} \right) + \varepsilon \nabla \cdot \left( \rho^{\varepsilon} h^{\varepsilon} {\rm {\bf v}}^{\varepsilon} \right) & = \nabla \cdot \left( {\rm k}^{\varepsilon} \nabla {\rm T}^{\varepsilon} \right) + \varepsilon ~ {\rm {\bf v}}^{\varepsilon} \cdot \nabla p^{\varepsilon} + \varepsilon^{3} ~ \boldsymbol{\tau}^{\varepsilon} : \nabla {\rm {\bf v}}^{\varepsilon} \\
 & + \varepsilon^{4} ~ {\rm {\bf v}}^{\varepsilon} \cdot \left( \nabla \cdot \left( \lambda \sigma \nabla \phi^{\varepsilon} \otimes \nabla \phi^{\varepsilon} \right) \right) \qquad \text{in } \Omega^{\varepsilon},
\end{aligned}
\end{equation}
\begin{equation}\label{ES4-6}
\partial_t \left( \rho_{S} {\rm c}_{p,S} {\rm T}_{S}^{\varepsilon} \right) = \nabla \cdot \left( {\rm k}_{S} \nabla {\rm T}_{S}^{\varepsilon} \right) \qquad \text{in } \Omega_{S}^{\varepsilon}.
\end{equation}
Along with the above dimensionless model equations, the corresponding boundary conditions on $\Gamma^{\varepsilon}$ are
\begin{equation}\label{ES4-7}
\nabla \phi^{\varepsilon} \cdot {\bf n}^{\varepsilon} = 0,
\end{equation}
\begin{equation}\label{ES4-8}
{\rm D}^{v} \rho^{\varepsilon} \nabla \chi^{v \varepsilon} \cdot {\bf n}^{\varepsilon} = 0,
\end{equation}
\begin{equation}\label{ES4-9}
{\rm {\bf v}}^{\varepsilon} = {\rm {\bf 0}},
\end{equation}
\begin{equation}\label{ES4-10}
{\rm T}^{\varepsilon} = {\rm T}_{S}^{\varepsilon},
\end{equation}
\begin{equation}\label{ES4-11}
{\rm k}^{\varepsilon} \nabla {\rm T}^{\varepsilon} \cdot {\bf n}^{\varepsilon} = {\rm k}_{S} \nabla {\rm T}_{S}^{\varepsilon} \cdot {\bf n}^{\varepsilon}.
\end{equation}

\subsection{Asymptotic Expansions}\label{Sec:6-3}
Using the homogenization ansatz as described in \hyperref[Sec:6-2]{Section 6.2}, we can write the unknowns as a series expansion in terms of the scale separation $\varepsilon$. For example, the phase field $\phi^{\varepsilon}$ can be written as
\begin{equation}\label{ES4-12}
\phi^{\varepsilon} (t, {\rm {\bf x}}) = \phi_0 (t, {\rm {\bf x}}, {\rm {\bf y}}) + \varepsilon \phi_1 (t, {\rm {\bf x}}, {\rm {\bf y}}) + \varepsilon^{2} \phi_2 (t, {\rm {\bf x}}, {\rm {\bf y}}) + \mathcal{O}(\varepsilon^{3}),
\end{equation}
where the functions $\phi_i (t, {\rm {\bf x}}, {\rm {\bf y}})$ are $Y$-periodic in ${\rm {\bf y}}$ and have explicit dependence on the micro- and the macro-scale variables. The micro-scale variable ${\rm {\bf y}}$ is defined locally in a pore $P$ and in solid $S$, while the macro-scale counterpart ${\rm {\bf x}}$ is in the entire (non-perforated) domain $\mathbb{D}$. Since $\phi^{\varepsilon}$ takes into account both the micro- and macro-scale behavior, we assume that the functions $\phi_i (t, {\rm {\bf x}}, {\rm {\bf y}})$ can separate between slow variability through ${\rm {\bf x}}$ and fast variability through ${\rm {\bf y}}$. Moreover, through $\phi_0$ one can capture the dominating behavior, while the subsequent terms describe the lower order behavior. We assume that the other dependent variables also have corresponding expansions.

Note that ${\rm {\bf y}}$ is a local variable, which is connected to the slow variable ${\rm {\bf x}}$ through the relation ${\rm {\bf y}} = \varepsilon^{-1} {\rm {\bf x}}$. As a consequence, the spatial derivatives need to be rewritten accordingly. Therefore, for a generic variable $\mathcal{V}$ one has
\begin{equation}\label{ES4-13}
\nabla \mathcal{V}({\rm {\bf x}}, {\rm {\bf y}}) = \nabla_{\rm {\bf x}} \mathcal{V}({\rm {\bf x}}, {\rm {\bf y}}) + \varepsilon^{-1} \nabla_{\rm {\bf y}} \mathcal{V}({\rm {\bf x}}, {\rm {\bf y}}),
\end{equation}
where $\nabla_{\rm {\bf x}}$ and $\nabla_{\rm {\bf y}}$ are the gradients w.r.t. ${\rm {\bf x}}$ and ${\rm {\bf y}}$, respectively. Now in order to find the dominating behavior, we insert these asymptotic expansions into the model equations \eqref{ES4-1}--\eqref{ES4-6} as well as in the boundary conditions \eqref{ES4-7})--\eqref{ES4-11}.

\subsubsection{Phase-field Equation}\label{Sec:6-3-1}
Equating the dominating $\mathcal{O}(1)$ terms from the phase-field equation \eqref{ES4-1}, gives
\begin{equation}\label{ES4-14}
\rho_0 \left\{ \partial_t \phi_0 + \nabla_{\rm {\bf y}} \cdot \left( \phi_0 {\rm {\bf v}}_0 \right) \right\} = \frac{\gamma}{\nu} \left\{ \nabla^{2}_{\rm {\bf y}} \phi_0 - \lambda^{-2} P'(\phi_0) \right\} - \sqrt{2} \lambda^{-1} \phi_0 (1 - \phi_0) f(\chi_{0}^{v}, {\rm T}_{0}) \qquad \text{in } P.
\end{equation}
The dominating term of the corresponding boundary condition \eqref{ES4-7} gives
\begin{equation}\label{ES4-15}
\nabla_{\rm {\bf y}} \phi_0 \cdot {\bf n}_P = 0 \qquad \text{on } \Gamma_P.
\end{equation}
Observe that the above equation is similar to the original phase-field equation \eqref{ES4-1}, but involves only spatial derivatives w.r.t. ${\rm {\bf y}}$. Although $\phi_0$ still depends on ${\rm {\bf x}}$, ${\rm {\bf x}}$ only appears as a parameter as no derivatives w.r.t. ${\rm {\bf x}}$ are involved. Recalling the $Y$-periodicity in ${\rm {\bf y}}$, $\phi_0$ solves the following \emph{cell problem} for the phase field
\begin{equation}\label{ES4-16}
\mathcal{P}^{1}_{\phi}: ~
\begin{cases}
    \rho_0 \left\{ \partial_t \phi_0 + \nabla_{\rm {\bf y}} \cdot \left( \phi_0 {\rm {\bf v}}_0 \right) \right\} = \frac{\gamma}{\nu} \left\{ \nabla^{2}_{\rm {\bf y}} \phi_0 - \lambda^{-2} P'(\phi_0) \right\} - \sqrt{2} \lambda^{-1} \phi_0 (1 - \phi_0) f(\chi_{0}^{v}, {\rm T}_0) ~~ \text{in } P, \\
    \nabla_{\rm {\bf y}} \phi_0 \cdot {\bf n}_P = 0 ~~ \text{on } \Gamma_P, \\
    \text{periodicity in } {\rm {\bf y}} ~ \text{across } \partial Y.
\end{cases}
\end{equation}
The cell problems are decoupled locally due to the periodicity requirement. That means, for every $\mathbf x\in\mathbb D$, we can solve a local cell problem with respect to $\mathbf y$, which is independent of the other cell problems.

\subsubsection{Mass Conservation Equations}\label{Sec:6-3-2}
The dominating $\mathcal{O}(1)$ terms in equation \eqref{ES4-2} give
\begin{equation}\label{ES4-17}
\partial_t \rho_0 + \nabla_{\rm {\bf y}} \cdot \left( \rho_0 {\rm {\bf v}}_0 \right) = 0 \qquad \text{in } P.
\end{equation}
When we integrate \eqref{ES4-17} with respect to ${\rm {\bf y}}$ over the domain $P$, apply Gauss' theorem in ${\rm {\bf y}}$, and use boundary conditions on $\Gamma_P$, we obtain
\begin{equation}\label{ES4-18}
\partial_t \overline{\rho_{0}} = 0 \qquad \text{in } \mathbb{D}.
\end{equation}
Here $\overline{()} := \frac{1}{|Y|} \int_P () ~d{\bf y}$ denotes the averaged value of the corresponding quantity. Next consider the $\mathcal{O}(\varepsilon)$ terms from equation \eqref{ES4-2}
\begin{equation}\label{ES4-19}
\partial_t \rho_1 + \nabla_{\rm {\bf x}} \cdot \left( \rho_0 {\rm {\bf v}}_0 \right) + \nabla_{\rm {\bf y}} \cdot \left( \rho_1 {\rm {\bf v}}_0 + \rho_0 {\rm {\bf v}}_1 \right) = 0 \qquad \text{in } P.
\end{equation}
Using similar arguments, integrating the above equation leads to
\begin{equation}\label{ES4-20}
\partial_t \overline{\rho_1} + \nabla_{\rm {\bf x}} \cdot \left( \overline{\rho_0 {\rm {\bf v}}_0} \right) = 0 \qquad \text{in } \mathbb{D}.
\end{equation}
Now, if we look at the $\mathcal{O}(1)+\mathcal{O}(\varepsilon)$ problem together, then we have
\begin{equation}\label{ES4-21}
\partial_t \overline{\rho_{\rm 0}} + \varepsilon \nabla_{\rm {\bf x}} \cdot \left( \overline{\rho_0 {\rm {\bf v}}_0} \right) = \mathcal{O}(\varepsilon) \qquad \text{in } \mathbb{D}.
\end{equation}
Note that as we are in the diffusion-dominated regime, the averaged leading order mass conservation equation \eqref{ES4-18} has no contribution from the convective term. One can certainly include the convective term in the mass conservation equation at the Darcy-scale, by considering the $\mathcal{O}(1) + \mathcal{O}(\varepsilon)$ problem, as done in equation \eqref{ES4-21}.

The dominating $\mathcal{O}(\varepsilon^{-2})$ term in equation \eqref{ES4-3} gives
\begin{equation}\label{ES4-22}
\nabla_{\rm {\bf y}} \cdot \left\{ {\rm D}^{v} \rho_{0} \nabla_{\rm {\bf y}} \chi_{0}^{v} \right\} = 0 \qquad \text{in } P.
\end{equation}
Integrating w.r.t ${\rm {\bf y}}$ over $P$, applying the Gauss theorem and the boundary condition \eqref{ES4-8} on $\Gamma_P$ together with periodicity, one gets
\begin{equation}\label{ES4-23}
\nabla_{\rm {\bf y}} \chi_{0}^{v} = 0 \qquad \text{in } P,
\end{equation}
meaning that $\chi_{0}^{v} = \chi_{0}^{v}(t, {\rm {\bf x}})$ is independent of ${\rm {\bf y}}$. The $\mathcal{O}(\varepsilon^{-1})$ terms in equation \eqref{ES4-3} read
\begin{equation}\label{ES4-24}
\nabla_{\rm {\bf y}} \cdot \left\{ {\rm D}^{v} \rho_{0} \left( \nabla_{\rm {\bf x}} \chi_{0}^{v} + \nabla_{\rm {\bf y}} \chi_{1}^{v} \right) \right\} = 0 \qquad \text{in } P.
\end{equation}
Now, due to linearity of the above equation, we can assume that
\begin{equation}\label{ES4-25}
\chi_{1}^{v} (t, {\rm {\bf x}}, {\rm {\bf y}}) = \widetilde{\chi}_{1}^{v} (t, {\rm {\bf x}}) + \sum_{j=1}^{d} \kappa^{j} (t, {\rm {\bf x}}, {\rm {\bf y}}) \partial_{x_{j}} \chi_{0}^{v} ~ .
\end{equation}
We note also that
\begin{equation}\label{ES4-26}
\nabla_{\rm {\bf x}} \chi_{0}^{v} = \sum_{j=1}^{d} e_{j} \partial_{x_{j}} \chi_{0}^{v} ~ .
\end{equation}
Using the above relations \eqref{ES4-25} and \eqref{ES4-26} the equation \eqref{ES4-24} can be rewritten, and the corresponding \emph{cell problem} becomes
\begin{equation}\label{ES4-27}
\mathcal{P}^{-1}_{\mathcal{D}}: \quad
\begin{cases}
\nabla_{\rm {\bf y}} \cdot \left\{ {\rm D}^{v} \rho_{0} (e_{j} + \nabla_{\rm {\bf y}} \kappa^{j}) \right\} = 0 \qquad \text{in } P,\\
{\rm D}^{v} \rho_{0} (e_{j} + \nabla_{\rm {\bf y}} \kappa^{j}) \cdot {\bf n}_{P} = 0 \qquad \text{on } \Gamma_P,\\
\text{periodicity in} ~ {\bf y} \qquad \text{across} ~ \partial Y.
\end{cases}
\end{equation}
Further the $\mathcal{O}(1)$ terms from equation \eqref{ES4-3} give
\begin{equation}\label{ES4-28}
\begin{aligned}
\partial_t \left( \rho_0 \chi_{0}^{v} \right) + \nabla_{\rm {\bf y}} \cdot \left( \rho_0 \chi_{0}^{v} {\rm {\bf v}}_{0} \right) = \nabla_{\rm {\bf x}} \cdot \left( {\rm D}^{v} \rho_{0} \nabla_{\rm {\bf x}} \chi_{0}^{v} \right) + \nabla_{\rm {\bf x}} \cdot \left( {\rm D}^{v} \rho_{0} \nabla_{\rm {\bf y}} \chi_{1}^{v} \right) + \nabla_{\bf y} \cdot \left( {\rm D}^{v} \rho_{0} \nabla_{\bf x} \chi_{1}^{v} \right) \\
+ \nabla_{\rm {\bf y}} \cdot \left( {\rm D}^{v} \rho_{1} \nabla_{\rm {\bf x}} \chi_{0}^{v} \right) + \nabla_{\rm {\bf y}} \cdot \left( {\rm D}^{v} \rho_{0} \nabla_{\rm {\bf y}} \chi_{2}^{v} \right) + \nabla_{\bf y} \cdot \left( {\rm D}^{v} \rho_{1} \nabla_{\bf y} \chi_{1}^{v} \right) \quad \text{in } P.
\end{aligned}
\end{equation}
To find an upscaled model we integrate equation \eqref{ES4-28} w.r.t ${ \rm {\bf y}}$ over the domain $P$, apply Gauss' theorem in ${\rm {\bf y}}$, and use the boundary condition on $\Gamma_{P}$. This leads to the upscaled mass conservation equation
\begin{equation}\label{ES4-29}
\partial_t \left( \overline{ \rho_{0} } \chi_{0}^{v} \right) = \nabla_{\rm {\bf x}} \cdot \left( \mathcal{D} \nabla_{\rm {\bf x}} \chi_{0}^{v} \right) \qquad \text{in } \mathbb{D},
\end{equation}
where the component $d_{ij} (t, {\rm {\bf x}})$ of the matrix $\mathcal{D}(t, {\rm {\bf x}})$ are given by
\begin{equation}\label{ES4-30}
d_{ij} (t, {\rm {\bf x}}) = \int_{P} {\rm D}^{v} \rho_{0} \left( \delta_{ij} + \partial_{y_i} \kappa^{j} \right) d{\rm {\bf y}}, \quad \text{with } i,j \in \{1, 2, \cdots, d \},
\end{equation}
where $\kappa^{j}$ is the solution of the cell problem \eqref{ES4-27}. Similarly, to include the convective term in the averaged mass conservation equation for vapor component, we again consider terms of $\mathcal{O}(\varepsilon)$, which would give
\begin{equation}\label{ES4-29EX}
\partial_t \left( \overline{ \rho_{0} } \chi_{0}^{v} \right) + \varepsilon \nabla_{\rm {\bf x}} \cdot \left( \overline{ \rho_0 {\rm {\bf v}}_{0} } \chi_0^{v}  \right) = \nabla_{\rm {\bf x}} \cdot \left( \mathcal{D} \nabla_{\rm {\bf x}} \chi_{0}^{v} \right) + \mathcal{O}(\varepsilon) \qquad \text{in } \mathbb{D}.
\end{equation}

\subsubsection{Momentum Balance Equation}\label{Sec:6-3-3}
The dominating $\mathcal{O}(\varepsilon^{-2})$ term in equation \eqref{ES4-4} yields
\begin{equation}\label{ES4-31}
\nabla_{\rm {\bf y}} p_0 = 0 \qquad \text{in } P,
\end{equation}
which implies that $p_0 = p_0 (t, {\rm {\bf x}})$ is independent of ${\rm {\bf y}}$. Next the $\mathcal{O}(\varepsilon^{-1})$ terms of equation \eqref{ES4-4} give
\begin{equation}\label{ES4-32}
\begin{aligned}
- \left( \nabla_{\rm {\bf x}} p_0 + \nabla_{\rm {\bf y}} p_1 \right) + \nabla_{\rm {\bf y}} \cdot \left\{ \mu_0 \left( \nabla_{\rm {\bf y}} {\rm {\bf v}}_0 + \nabla_{\rm {\bf y}} {\rm {\bf v}}_0^{T} \right) \right\} + \nabla_{\rm {\bf y}} \cdot \left\{ \xi_0 \left( \nabla_{\rm {\bf y}} \cdot {\rm {\bf v}}_0 \right) {\bf I} \right\} & \\
- \nabla_{\rm {\bf y}} \cdot \left( \lambda \sigma \nabla_{\rm {\bf y}} \phi_0 \otimes \nabla_{\rm {\bf y}} \phi_0 \right) + \rho_0 {\rm {\bf g}} & = 0 \qquad \text{in } P.
\end{aligned}
\end{equation}
We use the linearity of the variables $p_1$ and ${\rm {\bf v}}_0$ in the above equation and determine $p_1$ and ${\rm {\bf v}}_0$ in terms of (the gradient of) $p_0$. The following linear relationships then emerge
\begin{equation}\label{ES4-33}
p_1 (t, {\rm {\bf x}}, {\rm {\bf y}}) = \Pi^0 (t, {\rm {\bf x}}, {\rm {\bf y}}) + \sum_{j=1}^{d} \Pi^j (t, {\rm {\bf x}}, {\rm {\bf y}}) \partial_{x_j} p_0 (t, {\rm {\bf x}}),
\end{equation}
\begin{equation}\label{ES4-34}
{\rm {\bf v}}_0 (t, {\rm {\bf x}}, {\rm {\bf y}}) = - {\rm {\bf w}}^{0} (t, {\rm {\bf x}}, {\rm {\bf y}}) - \sum_{j=1}^{d} {\rm {\bf w}}^j (t, {\rm {\bf x}}, {\rm {\bf y}}) \partial_{x_j} p_0 (t, {\rm {\bf x}}).
\end{equation}
Now $\Pi^{j}(t, {\rm {\bf x}}, {\rm {\bf y}})$ and ${\rm {\bf w}}^{j}(t, {\rm {\bf x}}, {\rm {\bf y}})$ solve the \emph{cell problems}
\begin{equation}\label{ES4-35}
\mathcal{P}_{\mathcal{K}}^{-1}: ~~
\begin{cases}
\left( {\rm {\bf e}}_j + \nabla_{\rm {\bf y}} \Pi^j \right) + \nabla_{\rm {\bf y}} \cdot \left\{ \mu_0 \left( \nabla_{\rm {\bf y}} {\rm {\bf w}}^j + \nabla_{\rm {\bf y}} {\rm {\bf w}}^{j^{T}} \right) + \xi_0 \left( \nabla_{\rm {\bf y}} \cdot {\rm {\bf w}}^j \right) {\bf I} \right\} = 0 \quad \text{in } P, \\
\nabla_{\rm {\bf y}} \cdot \left( \rho_0 {\rm {\bf w}}^{j} \right) = 0 \quad \text{in } P, \\
{\rm {\bf w}}^{j} = {\rm {\bf 0}} \quad \text{on } \Gamma_P, \\
\text{periodicity in } {\rm {\bf y}} \text{ across } \partial Y \text{ and } \int_P \Pi^j d{\bf y} = 0, ~ j \in \left\{ 1, 2, \cdots, d \right\}.
\end{cases}
\end{equation}
In addition, we also have
\begin{equation}\label{ES4-36}
\mathcal{P}_{\mathcal{G}}^{-1}: ~~
\begin{cases}
\nabla_{\rm {\bf y}} \Pi^{0} = \nabla_{\rm {\bf y}} \cdot \left\{ \mu_0 \left( \nabla_{\rm {\bf y}} {\rm {\bf w}}^0 + \nabla_{\rm {\bf y}} {\rm {\bf w}}^{0^{T}} \right) +  \xi_0 \left( \nabla_{\rm {\bf y}} \cdot {\rm {\bf w}}^0 \right) {\bf I} \right\} \\
\qquad \quad - \nabla_{\rm {\bf y}} \cdot \left( \lambda \sigma \nabla_{\rm {\bf y}} \phi_0 \otimes \nabla_{\rm {\bf y}} \phi_0 \right) + \rho_0 {\rm {\bf g}} \quad \text{in } P, \\
\partial_t \rho_0 = \nabla_{\rm {\bf y}} \cdot \left( \rho_0 {\rm {\bf w}}^{0} \right) \quad \text{in } P,  \\
{\rm {\bf w}}^{0} = {\rm {\bf 0}} \quad \text{on } \Gamma_P. \\
\text{periodicity in } {\rm {\bf y}} \text{ across } \partial Y \text{ and } \int_P \Pi^0 d{\bf y} = 0.
\end{cases}
\end{equation}
Now multiplying the relation \eqref{ES4-34} by $\rho_{0}$ and then averaging over $Y$ leads to
\begin{equation}\label{ES4-37}
\overline{ \rho_{0} \rm {\bf v}}_{0} = - \mathcal{K} \nabla_{\rm {\bf x}} p_0 - \mathcal{G} \qquad \text{in } \mathbb{D},
\end{equation}
where the components of the \emph{permeability tensor} $\mathcal{K} (t, {\rm {\bf x}})$ are given by
\begin{equation}\label{ES4-38}
k_{ij} (t, {\rm {\bf x}}) = \int_P \rho_0 w_i^j d{\rm {\bf y}} \qquad \text{with } i,j \in \{1, 2, \cdots, d \},
\end{equation}
and the components of the column vector $\mathcal{G} (t, {\rm {\bf x}})$ are given by
\begin{equation}\label{ES4-39}
g_{i} (t, {\rm {\bf x}}) = \int_P \rho_0 w_i^0 d{\rm {\bf y}} \qquad \text{with } i \in \{1, 2, \cdots, d \}.
\end{equation}
Here, $w^j_i$ and $w^0_i$ are the components of ${\bf w}^j$ and ${\bf w}^0$, respectively, which are the solutions of the cell problems (\ref{ES4-35}) and (\ref{ES4-36}).

\subsubsection{Energy balance equation}\label{Sec:6-3-4}
The dominating $\mathcal{O}(\varepsilon^{-2})$ term from the energy balance equations \eqref{ES4-5} and \eqref{ES4-6}, respectively, gives
\begin{equation}\label{ES4-40}
\nabla_{\rm {\bf y}} \cdot \left( {\rm k}_0 \nabla_{\rm {\bf y}} {\rm T}_0 \right) = 0 \qquad \text{in } P,
\end{equation}
\begin{equation}\label{ES4-41}
\nabla_{\rm {\bf y}} \cdot \left( {\rm k}_S \nabla_{\rm {\bf y}} {\rm T}_{S0} \right) = 0 \qquad \text{in } S.
\end{equation}
The dominating $\mathcal{O}(1)$ and $\mathcal{O}(\varepsilon^{-1})$ terms from the corresponding boundary conditions \eqref{ES4-10} and \eqref{ES4-11}, respectively, lead to
\begin{equation}\label{ES4-42}
{\rm T}_0 = {\rm T}_{S0} \qquad \text{on } \Gamma_P,
\end{equation}
\begin{equation}\label{ES4-43}
{\rm k}_0 \nabla_{\rm {\bf y}} {\rm T}_0 \cdot {\bf n}_P + {\rm k}_S \nabla_{\rm {\bf y}} {\rm T}_{S0} \cdot {\bf n}_P = 0 \qquad \text{on } \Gamma_P,
\end{equation}
along with periodicity in ${\rm {\bf y}}$. This implies that ${\rm T}_0 = {\rm T}_{S0} = {\rm T}_0 (t, {\rm {\bf x}})$ is independent of ${\rm {\bf y}}$. Further, using this fact, the $\mathcal{O}(\varepsilon^{-1})$ terms from \eqref{ES4-5} and \eqref{ES4-6}, respectively, give
\begin{equation}\label{ES4-44}
\nabla_{\rm {\bf y}} \cdot \left\{ {\rm k}_0 \left( \nabla_{\rm {\bf x}} {\rm T}_0 + \nabla_{\rm {\bf y}} {\rm T}_1 \right) \right\} = 0 \qquad \text{in } P,
\end{equation}
\begin{equation}\label{ES4-45}
\nabla_{\rm {\bf y}} \cdot \left\{ {\rm k}_S \left( \nabla_{\rm {\bf x}} {\rm T}_{S0} + \nabla_{\rm {\bf y}} {\rm T}_{S1} \right) \right\} = 0 \qquad \text{in } S.
\end{equation}
As we have linearity in this equation, we can assume linear dependence of ${\rm T}_1$ and ${\rm T}_{S1}$ on the partial ${\rm {\bf x}}$-derivatives of ${\rm {\bf T}}_0$ in the following way
\begin{equation}\label{ES4-46}
{\rm T}_1 (t, {\rm {\bf x}}, {\rm {\bf y}}) = \widetilde{\rm T}_1 (t, {\rm {\bf x}}) + \sum_{j = 1}^{d} \psi^{j} (t, {\rm {\bf x}}, {\rm {\bf y}}) \partial_{x_{j}} {\rm T}_0 (t, {\rm {\bf x}}),
\end{equation}
\begin{equation}\label{ES4-47}
{\rm T}_{S1} (t, {\rm {\bf x}}, {\rm {\bf y}}) = \widetilde{\rm T}_{S1} (t, {\rm {\bf x}}) + \sum_{j = 1}^{d} \eta^{j} (t, {\rm {\bf x}}, {\rm {\bf y}}) \partial_{x_{j}} {\rm T}_0 (t, {\rm {\bf x}}),
\end{equation}
where the $\psi^{j}$ and $\eta^{j}$ are weight functions that need to be solved for in a specific cell ${\rm {\bf x}}$. Substituting ${\rm T}_1$ and ${\rm T}_{S1}$ in \eqref{ES4-44} and \eqref{ES4-45} with the new linear expressions \eqref{ES4-46} and \eqref{ES4-47} and requiring these equations to be fulfilled for any partial derivative of $\rm T_0$ yields the \emph{cell problem}
\begin{equation}\label{ES4-48}
\mathcal{P}^{-1}_{\mathcal{A}}: ~~
\begin{cases}
\nabla_{\rm {\bf y}} \cdot \left\{ {\rm k}_0 \left( e_j + \nabla_{\rm {\bf y}} \psi^{j} \right) \right\} = 0 \quad \text{in } P, \\
\nabla_{\rm {\bf y}} \cdot \left\{ {\rm k}_S \left( e_j + \nabla_{\rm {\bf y}} \eta^{j} \right) \right\} = 0 \quad \text{in } S,  \\ 
\psi^{j} = \eta^{j} \quad \text{on } \Gamma_{P}, \\
{\rm k}_0 \left( e_j + \nabla_{\rm {\bf y}} \psi^{j} \right) \cdot {\bf n}_P + {\rm k}_S \left( e_j + \nabla_{\rm {\bf y}} \eta^{j} \right) \cdot {\bf n}_P = 0 \quad \text{on } \Gamma_{P}, \\
\text{periodicity in } {\bf y} \quad \text{across } \partial Y.
\end{cases}
\end{equation}
Now $\mathcal{O}(1)$ terms from equations \eqref{ES4-5} and \eqref{ES4-6} lead to
\begin{equation}\label{ES4-49}
\begin{aligned}
\partial_t \left( \rho_0 u_0 \right) + \nabla_{\rm {\bf y}} \cdot \left( \rho_0 h_0 {\rm {\bf v}}_0 \right) = {\rm {\bf v}}_{0} \cdot \nabla_{\rm {\bf y}} p_{0} + \nabla_{\rm {\bf x}} \cdot \left( {\rm k}_0 \nabla_{\rm {\bf x}} {\rm T}_0 \right) + \nabla_{\rm {\bf x}} \cdot \left( {\rm k}_0 \nabla_{\rm {\bf y}} {\rm T}_1 + {\rm k}_1 \nabla_{\rm {\bf y}} {\rm T}_0 \right) \\
+ \nabla_{\rm {\bf y}} \cdot \left( {\rm k}_0 \nabla_{\rm {\bf x}} {\rm T}_1 + {\rm k}_1 \nabla_{\rm {\bf x}} {\rm T}_0 \right) + \nabla_{\rm {\bf y}} \cdot \left( {\rm k}_0 \nabla_{\rm {\bf y}} {\rm T}_2 + {\rm k}_1 \nabla_{\rm {\bf y}} {\rm T}_1 + {\rm k}_2 \nabla_{\rm {\bf y}} {\rm T}_0 \right) \qquad \text{in } P,
\end{aligned}
\end{equation}
\begin{equation}\label{ES4-50}
\begin{aligned}
\partial_t \left( \rho_S {\rm C}_{p,S} {\rm T}_{0} \right) = \nabla_{\rm {\bf x}} \cdot \left( {\rm k}_S \nabla_{\rm {\bf x}} {\rm T}_{0} \right) + \nabla_{\rm {\bf x}} \cdot \left( {\rm k}_S \nabla_{\rm {\bf y}} {\rm T}_{S1} \right) & + \nabla_{\rm {\bf y}} \cdot \left( {\rm k}_S \nabla_{\rm {\bf x}} {\rm T}_{S1} \right) \\
 & + \nabla_{\rm {\bf y}} \cdot \left( {\rm k}_S \nabla_{\rm {\bf y}} {\rm T}_{S2} \right) \qquad \text{in } S.
\end{aligned}
\end{equation}
Integrating the equations \eqref{ES4-49} and \eqref{ES4-50} w.r.t ${\rm {\bf y}}$ over $P$ and $S$, respectively, and then adding those leads to
\begin{equation}\label{ES4-51}
\partial_t \left( \overline{\rho_0} u_0 \right) + \partial_t \left( \rho_S {C}_{p,S} {\rm T}_{0} \right) = \nabla_{\rm {\bf x}} \cdot \left( \mathcal{A} \nabla_{\rm {\bf x}} {\rm T}_0 \right) \quad \text{in } \mathbb{D},
\end{equation}
where the component $a_{ij} (t, {\rm {\bf x}})$ of the matrix $\mathcal{A}(t, {\rm {\bf x}})$ are given by
\begin{equation}\label{ES4-52}
a_{ij} (t, {\rm {\bf x}}) = \int_{P} \left\{ {\rm k}_0 \left( \delta_{ij} + \partial_{y_i} \psi^{j} \right) + {\rm k}_S \left( \delta_{ij} + \partial_{y_i} \eta^{j} \right) \right\} d{\rm {\bf y}}, \quad \text{with } i,j \in \{1, 2, \cdots, d \},
\end{equation}
where $\psi^{j}$ and $\eta^{j}$ are the solution of the cell problem \eqref{ES4-48}. Again to include the convective term in the energy balance equation at the Darcy-scale, we consider terms of $\mathcal{O}(\varepsilon)$, which gives
\begin{equation}\label{ES4-51EX}
\partial_t \left( \overline{\rho_0} u_0 \right) + \partial_t \left( \rho_S {C}_{p,S} {\rm T}_{0} \right) + \varepsilon \nabla_{\rm {\bf x}} \cdot \left( \overline{ \rho_0 {\rm {\bf v}}_0 } h_0 \right) = \nabla_{\rm {\bf x}} \cdot \left( \mathcal{A} \nabla_{\rm {\bf x}} {\rm T}_0 \right) + \mathcal{O}(\varepsilon) \quad \text{in } \mathbb{D}.
\end{equation}
One may note that the internal energy ($u_0$) and the specific enthalpy ($h_0$) are functions of pressure ($p_0$) and temperature (${\rm T}_0$). Since, $p_0$ and ${\rm T}_0$ are not affected by the local variation, in other words, independent of ${\bf y}$, consequently, $u_0$ and $h_0$ are also independent of ${\bf y}$. Therefore, the averaging over ${\bf y}$ does not affect these quantities in equations \eqref{ES4-51} and \eqref{ES4-51EX}.

\section{Summary and Outlook}\label{Sec:7}
We start this section with summarizing the upscaled system of equations consisting of the conservation equations on the macro scale derived in \hyperref[Sec:6-3]{Section 6.3}. Further, the upscaled system is supplemented with the cell problems defined in the equations \eqref{ES4-27}, \eqref{ES4-35}, \eqref{ES4-36} and \eqref{ES4-48}, which are solved locally in each pore and that provide effective quantities for the upscaled model. When only accounting for the behavior of leading order $\mathcal{O}(1)$ on the macro scale ${\bf x} \in \mathbb{D}$, and for $t > 0$, the model equations are
\begin{equation*}
\begin{aligned}
\partial_t \overline{\rho}_{0} & = 0 & \quad \text{in } \mathbb{D},\\
\partial_t \left( \overline{ \rho_0 } \chi_{0}^{v} \right) & = \nabla_{\rm {\bf x}} \cdot \left( \mathcal{D} \nabla_{\rm {\bf x}} \chi_{0}^{v} \right) & \quad \text{in } \mathbb{D},\\
\partial_t \left( \overline{\rho_0} u_0 \right) + \partial_t \left( \rho_S {C}_{p,S} {\rm T}_{0} \right) & = \nabla_{\rm {\bf x}} \cdot \left( \mathcal{A} \nabla_{\rm {\bf x}} {\rm T}_0 \right) & \quad \text{in } \mathbb{D}.
\end{aligned}
\end{equation*}
Note that, since we are considering a diffusion-dominated process, the leading order upscaled equations have no contribution from the convective terms. However, these can be included in the upscaled equation through considering also $\mathcal{O}(\varepsilon)$ contributions at the macro scale ${\bf x} \in \mathbb{D}$, and for $t > 0$, which leads to
\begin{equation*}
\begin{aligned}
\partial_t \overline{\rho_{0}} + \varepsilon \nabla_{\rm {\bf x}} \cdot \left( \overline{\rho_0 {\rm {\bf v}}_0} \right) & = \mathcal{O}(\varepsilon) & \quad \text{in } \mathbb{D},\\
\partial_t \left( \overline{ \rho_{0} } \chi_{0}^{v} \right) + \varepsilon \nabla_{\rm {\bf x}} \cdot \left\{ \overline{ \rho_0 {\rm {\bf v}}_{0} } \chi_{0}^{v} \right\} & = \nabla_{\rm {\bf x}} \cdot \left( \mathcal{D} \nabla_{\rm {\bf x}} \chi_{0}^{v} \right) + \mathcal{O}(\varepsilon) & \quad \text{in } \mathbb{D},\\
\overline{\rho_{0} \rm {\bf v}}_{0} & = - \mathcal{K} \nabla_{\rm {\bf x}} p_0 - \mathcal{G} & \quad \text{in } \mathbb{D},\\
\partial_t \left( \overline{\rho_{0}} u_{0} \right) + \partial_t \left( \rho_S {C}_{p,S} {\rm T_{0}} \right) + \varepsilon \nabla_{\rm {\bf x}} \cdot \left( \overline{ \rho_0 {\rm {\bf v}}_0 } h_0 \right) & = \nabla_{\rm {\bf x}} \cdot \left( \mathcal{A} \nabla_{\rm {\bf x}} {\rm T}_0 \right) + \mathcal{O}(\varepsilon) & \quad \text{in } \mathbb{D}.
\end{aligned}
\end{equation*}
Here the phase field $\phi_0 (t, {\bf x}, {\bf y})$ is updated locally in each pore by solving the $\mathcal{P}^{1}_{\phi}$ problem given in equation \eqref{ES4-16} for all ${\bf x} \in \mathbb{D}$ and $t > 0$. Note that the velocity ${\rm{\bf v}}_0$ is obtained from the relation \eqref{ES4-34}. Further, the effective matrices $\mathcal{D}(t, {\bf x})$, $\mathcal{K} (t, {\bf x})$, $\mathcal{G} (t, {\bf x})$ and $\mathcal{A} (t, {\bf x})$ are computed through the expressions given in \eqref{ES4-30}, \eqref{ES4-38}, \eqref{ES4-39}, and \eqref{ES4-52}, along with the corresponding cell problems defined by $\mathcal{P}^{-1}_{\mathcal{D}}$, $\mathcal{P}^{-1}_{\mathcal{K}}$, $\mathcal{P}^{-1}_{\mathcal{G}}$, and $\mathcal{P}^{-1}_{\mathcal{A}}$ problems given in equations \eqref{ES4-27}, \eqref{ES4-35}, \eqref{ES4-36}, and \eqref{ES4-48}, respectively. Similar cell problems has already been implemented in various studies in the literature. We refer to Kelm \emph{et al.} \cite{Kelm2024} for the implementation of the similar flow cell problems, namely $\mathcal{P}^{-1}_{\mathcal{K}}$ and $\mathcal{P}^{-1}_{\mathcal{G}}$. Next, for the implementation of the similar $\mathcal{P}^{1}_{\phi}$ and $\mathcal{P}^{-1}_{\mathcal{D}}$ problems, we refer to Bringedal \emph{et al.} \cite{Bringedal2020}. Furthermore, similar cell problem related to energy balance, namely $\mathcal{P}^{-1}_{\mathcal{A}}$, can be found in Auriault \emph{et al.} \cite{Auriault2009}.

In summary, we have derived a phase-field model for evaporation in porous media, which explicitly accounts for the vapor component in the system. Moreover, we have used the most general form of the energy balance equation and considered heat transfer between the solid and the fluids. Compared to available phase-field models for evaporation and condensation, e.g.~\cite{Berti2009,Badillo2012,Dreyer2014,Fabrizio2016}, the present pore-scale phase-field model is more general and includes behavior relevant for evaporation in porous media. Unlike Badillo \cite{Badillo2012}, Berti and Giorgi \cite{Berti2009} and Dreyer \emph{et al.} \cite{Dreyer2014}, the present model explicitly considers the presence of a vapor component along with the liquid and gas phases in the system. The model equations consist of classical conservation equations of the continuum mechanics along with an additional balance equation for phase evolution/transition of Allen-Cahn type. Moreover, one may note that the energy balance equation is considered in the most general form. We have shown using matched asymptotic expansions that the phase-field model reduces to the expected sharp-interface model in the limit as the diffuse-interface width approaches zero. Consequently, the phase-field model captures mass, momentum and energy transport between the fluid phases during evaporation.

The major numerical challenge in the sharp-interface formulation is handling the time-dependent domains along with the moving interface. However, in the phase-field formulation the moving interface is approximated by a diffuse interface of very small width, which leads the model equations to be solved in a stationary domain. Hence, the proposed phase-field model simplifies the development of numerical simulation tools. The phase-field model developed here can be seen as a pore-scale model in case of a porous medium. Further, we have implemented the developed phase-field formulation in DuMux and studied numerical examples to investigate the efficiency of the formulation. We observed that the considered numerical examples are able to capture shrinking droplet and consequently the evaporation process. We noticed evaporation related cooling effect in both examples, which is evident from the steep temperature gradient from the droplet to the neighbouring phase. However, the pore-scale simulation is strongly influenced by the time-step size restriction due to the convergence issue with Newton's method for such a complex non-linear system. As this numerical challenge can be surmounted by rather implementing an upscaled model accounting for small-scale interface movements only, a Darcy-scale counterpart of the phase-field model is derived using homogenization techniques. We obtained cell problems providing effective mass diffusion, permeability, and heat conductivity. We believe that the proposed phase-field model and the corresponding Darcy-scale counterpart serve the long standing requirement of an effective model for porous-medium evaporation incorporating also pore-scale effects. Further, as the next step implementation of the resulting two-scale phase-field formulation is in development.  


\section*{Acknowledgment}
We thank the Deutsche Forschungsgemeinschaft (DFG, German Research Foundation) for supporting this work by funding SFB 1313, Project Number 327154368.



\end{document}